\documentclass[bj]{imsart}

\RequirePackage[OT1]{fontenc}
\RequirePackage{amsthm,amsmath,natbib}
\RequirePackage[colorlinks,citecolor=blue,urlcolor=blue]{hyperref}
\usepackage{graphicx, subfigure}
\usepackage{amssymb, dsfont, centernot}
\usepackage{enumerate,enumitem}
\usepackage{stmaryrd} 

\arxiv{arXiv:0000.0000}

\startlocaldefs
\numberwithin{equation}{section}
\theoremstyle{plain}
\newtheorem{thm}{Theorem}[section]
\theoremstyle{remark}

\newtheorem{rems}{Remarks}[section]
\endlocaldefs

\renewcommand{\P}{\mathbb P}
\newcommand{\R}{\mathbb R}
\newcommand{\E}{\mathbb E}
\newcommand{\indep}{\rotatebox[origin=c]{90}{$\models$}}
\newcommand{\nindep}{\centernot{\indep}}
\newcommand{\cov}{\text{\normalfont cov}}
\newcommand{\cor}{\text{\normalfont cor}}

\newcommand{\argmax}{\text{\normalfont argmax}}

\newtheorem{mydef}{Definition}[section]
\newtheorem{prop}{Proposition}[section]
\newtheorem{lemme}{Lemma}[section]

\begin{document}

\begin{center}
{\Large
	{\sc  Graph estimation for Gaussian data zero-inflated by double truncation.}
}
\bigskip

  Anne G\'egout-Petit$^{a}$, Aur\'elie Gueudin-Muller$^{a}$ \& Cl\'emence Karmann$^{a}$
\bigskip

{\it
$^{a}$ Universit\'e de Lorraine, CNRS, Inria, IECL, F-54000 Nancy, France, Inria BIGS Team ;  anne.gegout-petit@univ-lorraine.fr, aurelie.gueudin@univ-lorraine.fr, clemence.karmann@inria.fr
}
\end{center}
\bigskip

{\bf Abstract.} We consider the problem of graph estimation in a zero-inflated Gaussian model. In this model, zero-inflation is obtained by double truncation (right and left) of a Gaussian vector. The goal is to recover the latent graph structure of the Gaussian vector with observations of the zero-inflated truncated vector. We propose a two step estimation procedure. The first step consists in estimating each term of the covariance matrix by maximising the corresponding bivariate marginal log-likelihood of the truncated vector. The second one uses the graphical lasso procedure to estimate the precision matrix sparsity, which encodes the graph structure. We then state some theoretical convergence results about the convergence rate of the covariance matrix and precision matrix estimators. These results allow us to establish consistency of our procedure with respect to graph structure recovery. We also present some simulation studies to corroborate the efficiency of our procedure.
\smallskip


\section{Introduction}

Multivariate data analysis often involves describing and explaining the relationships among a set of
variables. Undirected graphical models offer a way to address this demand by using a graph to represent a model. A graph is a set of nodes and edges which can be represented as a graphic in order to make it easier to study, visually or computationally. Undirected graphical models are based on the conditional independence: a 
relation between two variables, represented by an edge in the graph, means that the corresponding variables are conditionally dependent given all the remaining variables. Among undirected graphical models, Gaussian graphical model provides a particularly convenient framework. This model assumes that the observations have a multivariate Gaussian distribution with mean $\mu$ and covariance matrix $\Sigma$. In this Gaussian setting, a direct relation between two variables corresponds to a non-zero entry in the precision matrix $\Sigma^{-1}$. In other words, if $\Sigma^{-1}_{ij}$ is zero, then variables $i$ and $j$ are conditionally independent given the other variables. Thus, graph estimation involves finding the zero pattern in the inverse covariance matrix. The theoretical approaches developed solve a maximum likelihood problem with an added $L_1$ penalty on the precision matrix to increase sparsity of the resulting graph. Many authors like \cite{YuanLin2007}, \cite{Dahl2008covariance} or \cite{banerjee2008model} used interior point methods to solve the exact maximisation of the $L_1$ penalised log-likelihood. One of the most powerful approach is the graphical lasso of \cite{friedman2008sparse}, who used a blockwise coordinate descent approach.\\

Furthermore, truncated Gaussian distributions received much attention in the second half of the last century. Cohen (\cite{Cohen49},  \cite{Cohen50}, \cite{Cohen57}) studied extensively mean and standard deviation estimation with likelihood maximisation for univariate doubly truncated Gaussian. These data correspond to Gaussian distributions which fall between two points of truncation $a$ and $b$, $a < b$. He distinguished cases when the number of ``unmeasured" observations, that is observations which fall into the tail(s), is known or not and whether we know the number in each tail. \cite{ShahJaiswal66} also studied this case by estimating parameters from first four sample moments.\\
Bivariate case was then naturally studied. \cite{DesRaj53} and \cite{Cohen55} studied mean, variance and covariance estimation when only one of the variables is truncated whereas \cite{Nath71}, \cite{Dyer73} then \cite{Muthen90} analysed it when both variables are truncated. In all these papers, as soon as one of the variables falls outside its points of truncation, none of the variables of the bivector is observed. In others words, data samples are only constituted of the Gaussian data for which the two variables are observed.\\
About multivariate case, \cite{Cohen57multi} estimate model parameters where only one variable is truncated by likelihood maximisation. \cite{Singh60} considers means and variances estimation in the case where only some variables are truncated. Later, \cite{GuptaTracy76}, \cite{Lee83} and \cite{ManjunathWilhelm09} studied moments when all the variables of the Gaussian vector are doubly truncated, that is, when all the variables of the Gaussian vector fall inside their points of truncation. Graph estimation and matrix covariance estimation do not seem to have yet been discussed in the literature.\\

In this paper, we address the problem of graph estimation in a zero-inflated Gaussian model. In this model, zero-inflation is obtained by double truncation (right and left) of a Gaussian vector. This means that each of the Gaussian variables are normally observed inside its points of truncation, but is null otherwise. If a variable is truncated, we then observe a zero instead of its value, but we still observe the other variables of the vector, contrary to the literature. Our goal is to recover the latent graph structure of the Gaussian vector, encoded in the precision matrix, with observations of the zero-inflated truncated vector. To retrieve this theoretical graph structure, we use the graphical lasso procedure which involves the empirical covariance matrix. Unlike the classic Gaussian setting, the Gaussian vector is not directly observed in our setting and its empirical covariance matrix is therefore inaccessible. We then propose another estimator for the covariance matrix whose theoretical guarantees, including the control of the convergence rate in infinite norm, required for the graphical lasso procedure, are studied.\\

The rest of the paper is organized as follows. In Section 2, we explicit the model and present the two step estimation procedure. The first step consists in estimating the covariance matrix, by estimating each term by maximising the corresponding bivariate marginal log-likelihood of the truncated vector, which is a non-convex optimisation problem. The second one relies on the graphical lasso to estimate the precision matrix. Section 3 contains two theoretical results about the convergence rate of the covariance and the precision matrix estimators. We first use recent results of \cite{Mei2017landscape} to set out that our covariance matrix estimator concentrates well in infinite norm around the theoretical covariance matrix. These results concern properties of the stationary points of non-convex empirical risk minimisation problems. Next, we use this first result and the consistency properties of graphical lasso, studied by \cite{Ravikumar2011} in a more general framework, to show consistency and sparsistency of our final estimator of the precision matrix. In Section 3, we also state the resulting theorem which establishes consistency of our procedure with respect to graph structure recovery. In Section 4, we present some simulations studies to corroborate its theoretical efficiency.

\section{Model and estimation procedure}

\subsection{Model}

Let $X$ be a Gaussian $p$-vector $X \sim \mathcal N_p(\mu, \Sigma^*)$ where $\mu = (\mu_j)_{j=1,\dots,p} \in \R^p$ is the mean vector and $\Sigma^*  = (\Sigma^*_{jk})_{1 \leq j,k \leq p} \in \mathcal M_p(\R)$ the covariance matrix. Let us consider the $p$-vector $Y$ defined as:  
\begin{align*}
Y_j = \mathds{1}_{a_j \leq X_j \leq b_j} X_j\text{ for all }j \in \{1, \dots,p\},
\end{align*}
where the points of truncation  $a_j, b_j \in \R$, $a_j < b_j$ are known and depend on $j$. The Gaussian vector $X$ is not directly observed, but it is observed through the zero-inflated truncated vector $Y$. Unlike what exists in the (multivariate) truncated Gaussian literature (for example, in \cite{GuptaTracy76}, \cite{Lee83} and \cite{ManjunathWilhelm09}), when one of the initial Gaussian variables falls outside its points of truncation, we observe a zero instead and we observe the rest of the vector according the same rule. In other words, our truncation does not consist in restricting to the observations of $X$ which falls into $[a_1,b_1] \times \dots \times [a_p,b_p] \subset \R^p$.\\

Assume that $\mu_j = 0$ and $\Sigma^*_{jj} = 1$ for all $j \in \{1, \dots, p\}$. In practice, we can estimate them with existing techniques for doubly truncated univariate Gaussian vector (for example, \citep{Cohen57}) if variables are not centered and scaled.

Gaussian graphical model is particularly appropriate for conditional dependency graph inference. Indeed, the precision matrix $\Theta^* := (\Sigma^*)^{-1}$ specifies the conditional dependency structure (see \cite{BookHastieTibshFriedman}). More precisely, the graph contains an edge between the variables $X_j$ and $X_k$ iff: 
\begin{equation}
\left. 
\begin{aligned}  \label{edges} X_j \longleftrightarrow X_k & \Longleftrightarrow  X_j\ \nindep\ X_k\ |\ (X_l)_{l \neq j,k}\\
 & \Longleftrightarrow  \cor(X_j ,X_k\ |\ (X_l)_{l \neq j,k}) \neq 0\\
 & \Longleftrightarrow  \Theta^*_{jk} \neq 0.
\end{aligned}
\right\}
\end{equation}

Given a symmetric positive definite matrix $M$, let us denote: 
\begin{equation} \label{edges}
E(M) = \{(j,k) \in \{1,\dots,p\}^2,\ j \neq k,\ M_{jk} \neq 0\}.
\end{equation} 
In particular, $E(\Theta^*)$ denotes the set of the edges of the theoretical graph.\\

The goal of this paper is to recover the latent graph structure of the variables of the Gaussian vector $X$ from observations of the zero-inflated truncated vector $Y$.

\subsection{Some theoretical tools}

To explicit the model and exhibit some complexities, we give some theoretical tools and will focus here on bivariate marginal likelihood from the truncated vector $Y$.

Let $(j,k) \in \{1,\dots,p\}^2,\ j<k$ and let $f_{jk}(x,y) = f(x,y,\Sigma^*_{jk})$ denotes the bivariate marginal log-likelihood function of the Gaussian couple $(X_j, X_k) \sim \mathcal N_2\Bigg(\begin{pmatrix} 
0 \\
0
\end{pmatrix}, \begin{pmatrix} 
1 & \Sigma^*_{jk} \\
\Sigma^*_{jk} & 1 
\end{pmatrix}\Bigg)$. 

With these notations, the likelihood of $(Y_j, Y_k)$ is then $\mathcal L_{jk}(\Sigma^*_{jk},y)$ where $y$ is an observation of the vector $Y$ and:
\begin{equation} \label{CoupleLikelihood}
\mathcal L_{jk}(\sigma,y) = \sum_{a, b = 0}^1 \phi_{a b, jk}(\sigma, y_j, y_k) n_{a b}(y_j, y_k), 
\end{equation}
with :
\begin{itemize}
\item  $n_{a b}(y_j, y_k) = \mathds{1}_{\zeta_j = a, \zeta_k = b}$ where $\zeta_l = \left\{
\begin{array}{l}
  1$ if $y_l \in [a_l,b_l]\setminus\{ 0\}, \\
  0$ if $y_l = 0.
\end{array}
\right.$
\item $\displaystyle \sum\limits_{a,b = 0}^1 n_{a b}(y_j, y_k) = 1$
\item $\phi_{11,jk}(\sigma, y_j, y_k) =  f(y_j,y_k,\sigma)$
\item $\phi_{01,jk}(\sigma, y_j, y_k) = \phi_{01,jk}(\sigma,y_k) = \displaystyle \int_{[a_j,b_j]^c} f(x,y_k,\sigma) \mathrm{d}x $
\item $\phi_{10,jk}(\sigma, y_j, y_k) = \phi_{10,jk}(\sigma,y_j) = \displaystyle \int_{[a_k,b_k]^c} f(y_j,y,\sigma) \mathrm{d}y $
\item $\phi_{00,jk}(\sigma, y_j, y_k) = \phi_{00,jk}(\sigma) = \displaystyle \iint_{[a_j,b_j]^c\times[a_k,b_k]^c} f(x,y,\sigma) \mathrm{d}x \mathrm{d}y $.\\
\end{itemize}

The likelihood (and log-likelihood) of a couple of variables of $Y$ involves four terms according to the nullity of each of the components of the couple. In the same way, the likelihood of $Y$ would involve $2 ^p$ terms by distinguishing all possible cases: the density of the Gaussian vector $X$ (no component of $Y$ is null), $p$ simple integrals (only one null component), $\dbinom{p}{2}$ double integrals (two null components), $\dbinom{p}{3}$ triple integrals, ..., one $p$-multiple integral (all the components are null). Writing the likelihood of the vector $Y$ becomes than quite complicated and tedious. This is why we choose to restrict to the study of the likelihoods of couples for the estimation.\\

In practice, we have a $n$-sample $\textbf{Y} := (Y^{(1)}, \dots, Y^{(n)})$ of the vector $Y$. The likelihood of the $n$-sample $\big((Y^{(i)}_j, Y^{(i)}_k)\big)_{i = 1, \dots,n}$ is then $\mathcal L^{(n)}_{jk}(\Sigma^*_{jk},\textbf{y})$ defined by:
\begin{align*}
\mathcal L^{(n)}_{jk}(\sigma,\textbf{y}) & = \prod_{i = 1}^n \mathcal L_{jk}(\sigma,y^{(i)}),\\
& = \prod_{i = 1}^n \sum_{a, b = 0}^1 \phi_{a b, jk}(\sigma, y_j^{(i)}, y_k^{(i)}) n_{a b}(y_j^{(i)}, y_k^{(i)}),
\end{align*}
where $\textbf{y} := (y^{(1)}, \dots, y^{(n)})$ is the realisation (value) of the $n$-sample $\textbf{Y}$. The log-likelihood is then $L^{(n)}_{jk}(\Sigma^*_{jk},\textbf{y})$ where:
\begin{align*}
L^{(n)}_{jk}(\sigma,\textbf{y})& =  \sum_{i = 1}^n \sum_{a, b = 0}^1 n_{a b}(y_j^{(i)}, y_k^{(i)}) \log\Big(\phi_{a b, jk}(\sigma, y_j^{(i)}, y_k^{(i)})\Big) \\
& =  \sum\limits_{\substack{i=1 \\ i : y_j^{(i)} = y_k^{(i)} = 0}}^{n}{\log\Big(\phi_{00, jk}(\sigma)\Big)} + \sum\limits_{\substack{i=1 \\ i : y_j^{(i)} = 0, y_k^{(i)} \neq 0}}^{n}{\log\Big(\phi_{01, jk}(\sigma,y_k^{(i)})\Big)} \\ 
& + \sum\limits_{\substack{i=1 \\ i : y_j^{(i)} \neq 0, y_k^{(i)} = 0}}^{n}{\log\Big(\phi_{10, jk}(\sigma,y_j^{(i)})\Big)} + \sum\limits_{\substack{i=1 \\ i : y_j^{(i)} \neq 0, y_k^{(i)} \neq 0}}^{n}{\log\Big(\phi_{11, jk}(\sigma,y_j^{(i)},y_k^{(i)})\Big)}. 
\end{align*}

\subsection{Estimation procedure} \label{ssect:EstiMat}

Our goal is to recover the latent graph structure of the Gaussian vector $X$, encoded in the precision matrix $\Theta^*$, from observations of the truncated vector $Y$. Our estimation procedure is a two step procedure. In the first instance, we estimate the covariance matrix of the Gaussian vector $X$. Then, we estimate the precision matrix by using the graphical lasso procedure \citep{friedman2008sparse} to recover the underlying graph structure.

\subsubsection{Step 1: covariance matrix estimation}

Estimating the covariance matrix $\Sigma^*$ of $X$ as the empirical covariance matrix of the $n$-sample $\textbf{Y}$ would lead to poor results because of the zero-inflation. \\

Another idea could be to maximise the likelihood of the vector $Y$. But we have seen that this likelihood involves $2^p$ terms and is too tedious.

Because of these difficulties, we estimate the covariance matrix by estimating each of its entries separately using the likelihood of the couples $(Y_j, Y_k)$, $j<k$ defined in \eqref{CoupleLikelihood}. 
More precisely, we estimate $\Sigma^*$ by $\widetilde \Sigma^{(n)}$ by estimating each of its entries $\Sigma^*_{jk}$ by maximisation of the log-likelihood of the $n$-sample $\big((Y^{(i)}_j, Y^{(i)}_k)\big)_ {i = 1, \dots, n}$ of the couple $(Y_j, Y_k)$, which is not convex.

\begin{mydef}[Estimator $\widetilde \Sigma^{(n)}$ of $\Sigma^*$]
The estimator $\widetilde \Sigma^{(n)} = (\widetilde \Sigma^{(n)}_{jk})_{1\leq j,k\leq p}$ of the covariance matrix $\Sigma^*$ is defined by: 
\begin{align}
\widetilde \Sigma^{(n)}_{jk} & = \underset{|\sigma| \leq 1}{\argmax}\ L^{(n)}_{jk}(\sigma,\textbf{y})\label{EstiSigma}\\
& = \underset{|\sigma| \leq 1}{\argmax}\ \dfrac{1}{n}L^{(n)}_{jk}(\sigma,\textbf{y}),\nonumber
\end{align} 
for all $j<k$, where $\textbf{y} := (y^{(1)}, \dots, y^{(n)})$ is the realisation of the $n$-sample $\textbf{Y}$.
\end{mydef}

\subsubsection{Step 2: precision matrix estimation}

As our goal is to recover the conditional dependency graph, it is natural to use the estimator of the precision matrix $\Theta^*$ given by the graphical lasso \citep{friedman2008sparse}. The graphical lasso is a procedure used in the Gaussian graphical model. It consists in estimating the precision matrix by maximising the penalised log-likelihood of the Gaussian model over the set $p \times p$ non-negative definite symmetric matrices: 
$$\underset{\Theta \succ 0}{\argmax}\ \log \det(\Theta)-\text{trace}(\Theta S) - \lambda_n ||\Theta||_{1,\text{off}},$$ 
where $||\Theta||_{1,\text{off}} = \sum\limits_{\substack{j,k =1 \\ j \neq k}}^p |\Theta_{jk}|$, $S$ is the empirical covariance matrix of $X$ and $\lambda_n > 0$ the regularisation parameter. This optimisation problem is convex and has an unique solution \citep{Ravikumar2011}.\\ 

In our case, the empirical covariance matrix of $X$ is not directly available. Instead of obtaining this matrix as the empirical covariance matrix of $Y$, we replace the empirical covariance matrix $S$ by the estimator $\widetilde \Sigma^{(n)}$ of $\Sigma^*$ obtained at the step 1 \eqref{EstiSigma}.

\begin{mydef}[Estimator $\widetilde \Theta^{(n)}$ of $\Theta^*$]
The estimator $\widetilde \Theta^{(n)}$ of the precision matrix $\Theta^*$ is defined as the unique solution of the following convex optimisation problem:
\begin{equation}\label{LassoRav}
\widetilde \Theta^{(n)} = \underset{\Theta \succ 0}{\argmax}\ \log \det(\Theta)-\text{trace}(\Theta\widetilde \Sigma^{(n)}) - \lambda_n ||\Theta||_{1,\text{off}}.
\end{equation}
\end{mydef}

Theoretical results of Section \ref{sect:ResultatsTheo} relate the estimators $\widetilde \Sigma^{(n)}$ and $\widetilde \Theta^{(n)}$ respectively defined in \eqref{EstiSigma} and \eqref{LassoRav} when the points of truncation $(a_j)_{1 \leq j \leq p}$ and $(b_j)_{1 \leq j \leq p}$ are known. 

\section{Convergence results} \label{sect:ResultatsTheo}

The goal of this Section is to show that the estimation procedure proposed in Subsection \ref{ssect:EstiMat} has strong theoretical guarantees. For that, we study theoretical properties of the estimator $\widetilde \Theta^{(n)}$ with regard to the recovery of the graph structure. Assume that the points of truncation $(a_j)_{1 \leq j \leq p}$ and $(b_j)_{1 \leq j \leq p}$ are known. 

\subsection{Covariance matrix estimator}

\subsubsection{Convergence rates in elementwise infinite norm}

In a first place, we give a result about the estimator $\widetilde \Sigma^{(n)}$ of the covariance matrix $\Sigma^*$ given by \eqref{EstiSigma}. Let us first set out two assumptions: 

\newcounter{saveenum}

\begin{enumerate}[label = \textbf{(H\arabic*)}, leftmargin = *]
  \item \label{majoration} For all $j < k$, $|\Sigma^*_{jk}| \neq 1$. Thus, there exists $\delta > 0$ such that for all $j < k$, $|\Sigma^*_{jk}| < 1 - \delta$. \\
  
  \item \label{Morse} Let $j < k$ and consider the application $g : \sigma \in [-1+\delta, 1-\delta] \mapsto \E\Big(L^{(n)}_{jk}(\sigma, \textbf{y})\Big)$. Then, we assume that:
\begin{itemize}[label = $\bullet$]
\item $-1+\delta$ and $1-\delta$ are not critical points of $g$,
\item $g$ has a finite number of critical points,
\item every critical points of $g$, different from $\Sigma^* _{jk}$, are non-degenerate, i.e.:
$$\text{for all $\sigma \neq \Sigma^*_{jk}$, }  g'(\sigma) = 0 \Rightarrow g''(\sigma) \neq 0.$$
\end{itemize}
Note that $\Sigma^*_{jk}$ is a non-degenerate critical point of $g$. This will be proved in the proof of Proposition \ref{Prop1} (see equations \eqref{ptcritique}).\\
  \setcounter{saveenum}{\value{enumi}}
\end{enumerate}

Proposition \ref{Prop1} states rate convergence results about the estimator $\widetilde \Sigma^{(n)}$ of the covariance matrix $\Sigma^*$ by bounding $\big|\big|\widetilde \Sigma^{(n)} - \Sigma^*\big|\big|_\infty$ with high probability. 

\begin{prop} \label{Prop1} Assume \ref{majoration} and \ref{Morse} and let $0 < \rho < 1$. There exist some known constants $B$, $C$ and $D$ such that letting $\dfrac{n}{\log n} \geq C \log\Big(\dfrac{B}{\rho} \Big)$, then the estimator of the covariance matrix $\widetilde \Sigma^{(n)}$ defined by \eqref{EstiSigma} satisfies:
\begin{align*}
\P\Bigg( \big|\big|\widetilde \Sigma^{(n)} - \Sigma^*\big|\big|_\infty  \geq D \sqrt{\dfrac{\log n}{n} \log\Big(\dfrac{B}{\rho} \Big)}  \Bigg) \leq \dfrac{p(p-1)}{2}\rho,\\
\end{align*}
where $||A||_\infty = \underset{j,k \in \{1,\dots,p\}}{\max} |A_{jk}|$ is the elementwise infinite norm of the matrix $A$. 
\end{prop}

\subsubsection{Proof of Proposition \ref{Prop1}}

Proof relies on Theorem 2 of \cite{Mei2017landscape}, who study the properties of the stationary points of non-convex empirical risk minimisation problems. We begin with three auxiliary Lemmas, proved in Appendix, which all state properties about the bivariate marginal likelihood defined in \eqref{CoupleLikelihood} or components of it: 

\begin{lemme} \label{L1} There exists $\gamma > 0$ such that, for all $j < k$, if $(y_j, y_k) \in [a_j, b_j]\times[a_k, b_k]$ and $\sigma \in [-1+\delta, 1-\delta]$, then for all $a, b \in \{0,1\}$, $\phi_{a b, jk}(\sigma, y_j, y_k) \geq \dfrac{1}{\gamma}$.
\end{lemme}

\begin{lemme} \label{L2} There exist $L_1, L_2 \text{ and } L_3 > 0$ such that for all $j < k$, if $(y_j, y_k) \in [a_j, b_j]\times[a_k, b_k]$ and $\sigma \in [-1+\delta, 1-\delta]$, then for all $a, b \in \{0,1\}$,
$$ \Big|\partial_\sigma^m \phi_{a b, jk}(\sigma, y_j, y_k)\Big| \leq L_m \text{, for $m \in \{1,2,3\}$.}  $$
\end{lemme}

\begin{lemme} \label{L3}  Let $j < k$.
\vspace*{0.01cm}
\begin{enumerate}
\item \label{L3.1} For all $\sigma \in [-1+\delta, 1-\delta]$ and for all $l \in \mathbb N^*$, 
$$\int_{\R^2} \partial_{\sigma}^l \mathcal L_{jk}(\sigma, y) \mathrm{d}\mu(y) = \partial_{\sigma}^l \displaystyle \int_{\R^2} \mathcal L_{jk}(\sigma, y) \mathrm{d}\mu(y) = 0,$$
 where $\mu$ is the measure on $\R^2$ defined by:
\begin{equation}\label{mu}
\mu := \delta_{0}\otimes\delta_0 + \delta_{0}\otimes\lambda + \lambda\otimes\delta_{0} + \lambda\otimes\lambda,
\end{equation}
where $\delta_a$ denotes the Dirac measure in $a \in \R$ and $\lambda$ the Lebesgue measure on $\R$.
\item \label{L3.2} For all $\sigma \in [-1+\delta, 1-\delta]$ and for all $l \in \mathbb N^*$,
 $$\partial_{\sigma}^l \E_{\Sigma^*_{jk}}\Big(\log \mathcal L_{jk}(\sigma, Y)\Big)  = \E_{\Sigma^*_{jk}}\Big(\partial_{\sigma}^l \log \mathcal L_{jk}(\sigma, Y)\Big),$$
 i.e., 
 $$\partial_{\sigma}^l \int_{\R^2} \log \mathcal L_{jk}(\sigma, y) \mathcal L_{jk}(\Sigma^*_{jk}, y) \mathrm{d}\mu(y)  = \int_{\R^2} \partial_{\sigma}^l \Big(\log \mathcal L_{jk}(\sigma, y)\Big)  \mathcal L_{jk}(\Sigma^*_{jk}, y) \mathrm{d}\mu(y).$$
\end{enumerate}
\end{lemme}


Fix $j<k$. With notations of \cite{Mei2017landscape}, let us set: 
\begin{align}
\ell_{jk}(\sigma, \textbf{y}) & =  \log \mathcal L_{jk}(\sigma, \textbf{y})  \nonumber \\
& = \sum_{a = 0}^1 \sum_{b = 0}^1 n_{a b}(y_j, y_k) \log\Big(\phi_{a b, jk}(\sigma, y_j, y_k)\Big) \label{ell}\\
\hat R_n(\sigma, \textbf{Y}) & = \dfrac{1}{n}L^{(n)}_{jk}(\sigma, \textbf{Y}) = \dfrac{1}{n} \sum_{i=1}^n \ell(\sigma, Y^{(i)})\\
R(\sigma) & = \E_{\Sigma^*_{jk}}\Big(\hat R_n(\sigma, \textbf{Y})\Big) = \E_{\Sigma^*_{jk}}\Big(\ell(\sigma, Y)\Big).\label{R}
\end{align} 

\begin{rems} \begin{itemize}[leftmargin = *]
\item To lighten notations, we drop the underscripts $jk$ and simply write $\ell$, $\hat R_n$ and $R$ instead of $\ell_{jk}$, $\hat R_{n, jk}$ and $R_{jk}$. 
\item Point \ref{L3.2} of Lemma \ref{L3} can be rewritten as:  \\
For all $\sigma \in [-1+\delta, 1-\delta]$ and for $l \in \mathbb N^*$:
\begin{equation} \label{L3bis}
R^{(l)}(\sigma) = \partial_{\sigma}^l \E_{\Sigma^*_{jk}}\Big(\ell(\sigma, Y))\Big)  = \E_{\Sigma^*_{jk}}\Big(\partial_{\sigma}^l \ell(\sigma, Y))\Big).
\end{equation}
\end{itemize}
\end{rems}

Theorem 2 of \cite{Mei2017landscape} requires four assumptions. Let us check these assumptions.\\

\textbf{(i) Gradient statistical noise.} \textit{The gradient of $\ell$ w.r.t. $\sigma$ is $\tau^2$-sub-Gaussian for some $\tau > 0$, i.e.: 
$$ \forall \sigma \in [-1+\delta, 1-\delta],\ \forall \lambda \in \R,\ \E\Bigg[\exp\bigg(\lambda\Big(\partial_\sigma \ell(\sigma,Y) - \E\big(\partial_\sigma \ell(\sigma, Y)\big)\Big) \bigg) \Bigg] \leq \exp\Big(\dfrac{\tau^2\lambda^2}{2}\Big). $$
}

Indeed, for all $y \in \prod\limits_{j=1}^p [a_j,b_j]$ and $\sigma \in [-1+\delta, 1-\delta]$,
\begin{align*}
\partial_\sigma \ell(\sigma, y) & =  \sum_{a = 0}^1 \sum_{b = 0}^1 n_{a b}(y_j, y_k) \dfrac{\partial_\sigma \phi_{a b, jk}(\sigma, y_j, y_k)}{\phi_{a b, jk}(\sigma, y_j, y_k)} \\ 
\text{Thus: } \Big| \partial_\sigma \ell(\sigma, y) \Big| & \leq \sum_{a = 0}^1 \sum_{b = 0}^1 n_{a b}(y_j, y_k) \dfrac{\big|\partial_\sigma \phi_{a b, jk}(\sigma, y_j, y_k)\big|}{\big|\phi_{a b, jk}(\sigma, y_j, y_k)\big|}\\
& \leq \sum_{a = 0}^1 \sum_{b = 0}^1 n_{a b}(y_j, y_k) \gamma L_1 = \gamma L_1 \text{ by Lemmas \ref{L1} and \ref{L2}.}
\end{align*}

This way, $\partial_\sigma \ell(\sigma,Y) - \E(\partial_\sigma \ell(\sigma, Y))$ is zero-mean and bounded by $2\gamma L_1$. By Theorem 9.9 of \cite{Stromberg}, $\partial_\sigma \ell(\sigma,Y) - \E(\partial_\sigma \ell(\sigma, Y))$ is then $\tau^2$-sub-Gaussian for $\tau = 2 \gamma L_1$. Assumption ``Gradient statistical noise" is satisfied.\\

\vspace{0.7cm}

\textbf{(ii) Hessian statistical noise.} \textit{The second derivative of $\ell$ w.r.t. $\sigma$ is $\tau^2$-sub-exponential ($\tau = 2\gamma L_1$), that is: 
$$|| \partial^2_\sigma \ell(\sigma,Y) ||_{\psi_1} \leq \tau^2 ,$$
where $||.||_{\psi_1}$ is the Orlicz $\psi_1$-norm defined by $||X||_{\psi_1} := \underset{k \geq 1}{\sup}\ \dfrac{1}{k} \E\Big(\big| X - \E(X) \big|^k \Big)^{\frac{1}{k}}$.\\
}

For all $y \in \prod\limits_{j=1}^p [a_j,b_j]$ and $\sigma \in [-1+\delta, 1-\delta]$,
\begin{align}
\partial^2_\sigma \ell(\sigma, y) & = \sum_{a = 0}^1 \sum_{b = 0}^1 n_{a b}(y_j, y_k) \bigg(\dfrac{\partial^2_\sigma \phi_{a b, jk}(\sigma, y_j, y_k)}{\phi_{a b, jk}(\sigma, y_j, y_k)} - \Big(\dfrac{\partial_\sigma \phi_{a b, jk}(\sigma, y_j, y_k)}{\phi_{a b, jk}(\sigma, y_j, y_k)}\Big)^2\bigg) \nonumber\\ 
\text{Thus: } \Big| \partial^2_\sigma \ell(\sigma, y) \Big| & \leq \sum_{a = 0}^1 \sum_{b = 0}^1 n_{a b}(y_j, y_k) \bigg(\dfrac{\big|\partial^2_\sigma \phi_{a b, jk}(\sigma, y_j, y_k)\big|}{\big|\phi_{a b, jk}(\sigma, y_j, y_k)\big|} + \Big(\dfrac{\big|\partial_\sigma \phi_{a b, jk}(\sigma, y_j, y_k)\big|}{\big|\phi_{a b, jk}(\sigma, y_j, y_k)\big|}\Big)^2\bigg)\nonumber\\
& \leq \sum_{a = 0}^1 \sum_{b = 0}^1 n_{a b}(y_j, y_k) (\gamma L_2 + \gamma^2L_1^2) \text{ by Lemmas \ref{L1} and \ref{L2}}\nonumber\\
& = \gamma L_2 + \gamma^2L_1^2. \label{ii}
\end{align}

Therefore, $\partial^2_\sigma \ell(\sigma,Y) - \E(\partial^2_\sigma \ell(\sigma, Y))$ is bounded by $2(\gamma L_2 + \gamma^2L_1^2)$ and for all $k \geq 1$, 
$$\dfrac{1}{k} \E\Big(\big| \partial^2_\sigma \ell(\sigma,Y) - \E(\partial^2_\sigma \ell(\sigma,Y)) \big|^k \Big)^{\frac{1}{k}} \leq \dfrac{2}{k}(\gamma L_2 + \gamma^2L_1^2) .$$
Hence, $||\partial^2_\sigma \ell(\sigma,Y) ||_{\psi_1} \leq 2(\gamma L_2 + \gamma^2L_1^2) \leq \tau^2 = 4\gamma^2L_1^2$ (we can possibly choose $L_1$ and $\gamma$ larger). So, $\partial^2_\sigma \ell(\sigma,Y)$ is $\tau^2$-sub-exponential. Assumption ``Hessian statistical noise" is satisfied.

\vspace{0.7cm}

\textbf{(iii) Hessian regularity.} \textit{\begin{enumerate}
\item The second derivative of $R$ (defined in \eqref{R}) is bounded at one point:
$$\text{there exists } |\sigma^*| \leq 1 - \delta \text{ and } H > 0 \text{ such that } \Big|R''(\sigma^*) \Big| \leq H.$$
\item The second derivative of $\ell$ w.r.t. $\sigma$ is Lipschitz continuous with integrable Lipschitz constant (w.r.t. $y$), i.e.:
$$\text{there exists } J^* > 0 \text{ such that } \E[J(Y)] \leq J^*,$$
where $J(y) = \underset{\substack{|\sigma_1|,|\sigma_2| \leq 1 - \delta \\ \sigma_1 \neq \sigma_2}}{\sup} \dfrac{|\partial^2_\sigma\ell(\sigma_1, y) - \partial^2_\sigma\ell(\sigma_2, y)|}{|\sigma_1 - \sigma_2|}$.
\item Constants $H$ and $J^*$ satisfy: $H \leq \tau^2$ and $J^* \leq \tau^3$. \\
\end{enumerate}
}

First, $R''(\sigma)  = \E_{\Sigma^*_{jk}}\Big(\partial_{\sigma}^2 \ell(\sigma, Y)\Big)$ by the point \ref{L3.2} of Lemma \ref{L3} and \eqref{L3bis}. By \eqref{ii}, $\Big| \partial^2_\sigma \ell(\sigma, Y) \Big| \leq  \gamma L_2 + \gamma^2L_1^2$ for all $\sigma \in [-1+\delta, 1-\delta]$, thus any $|\sigma^*| \leq 1 - \delta$ and $H = \gamma L_2 + \gamma^2L_1^2$ are appropriate. Moreover, we have $H \leq \tau^2 = 4\gamma^2L_1^2$ (with $L_1$ and $\gamma$ possibly chosen larger).\\

For all $y \in \prod\limits_{j=1}^p [a_j,b_j]$ and $\sigma \in [-1+\delta, 1-\delta]$, we have (with a slight lightening of notations): 
\begin{align*}
\partial^3_\sigma \ell(\sigma, y) & = \sum_{a = 0}^1 \sum_{b = 0}^1 n_{a b}(y_j, y_k) \bigg(\dfrac{\partial^3_\sigma \phi_{a b, jk}}{\phi_{a b, jk}} - 3\dfrac{\partial_\sigma \phi_{a b, jk}\partial^2_\sigma \phi_{a b, jk}}{\phi^2_{a b, jk}} + 2\Big(\dfrac{\partial_\sigma \phi_{a b, jk}}{\phi_{a b, jk}}\Big)^3\bigg) \\ 
\text{Thus: } \Big| \partial^3_\sigma \ell(\sigma, y) \Big| & \leq \sum_{a = 0}^1 \sum_{b = 0}^1 n_{a b}(y_j, y_k) (\gamma L_3 + 3\gamma^2L_1L_2 + 2\gamma^3L_1^3) \text{ (Lemmas \ref{L1} and \ref{L2})}\\
& = \gamma L_3 + 3\gamma^2L_1L_2^2 + 2\gamma^3L_1^3.
\end{align*}

Therefore, for all $y \in \prod\limits_{j=1}^p [a_j,b_j]$, $J(y) \leq \gamma L_3 + 3\gamma^2L_1L_2^2 + 2\gamma^3L_1^3$ by the mean value theorem. Taking $J^* = \gamma L_3 + 3\gamma^2L_1L_2^2 + 2\gamma^3L_1^3$ with $L_1$ and $\gamma$ possibly chosen larger, we have $J^* \leq \tau^3 = 8\gamma^3L_1^3$. Assumption ``Hessian regularity" is satisfied.

\vspace{0.7cm}

\textbf{(iv) Morse.} \textit{There exists $\epsilon > 0$ and $\eta > 0$ such that $R$ is $(\epsilon,\eta)$ strongly Morse, i.e.:
\begin{enumerate}
\item  $|R'(\sigma)| > \epsilon$ for all $\sigma$ such that $|\sigma| = 1-\delta$ and,
\item  for all $\sigma$ such that $|\sigma| < 1 - \delta$: 
$$ |R'(\sigma)| \leq \epsilon \Rightarrow |R''(\sigma)| \geq \eta.$$
\end{enumerate}
}

In other words, $R$ satisfies this assumption if $-1+\delta$ and $1-\delta$ are not critical points of $R$ and if $R$ has a finite number of critical points, which are moreover non-degenerate:
$$ R'(\sigma) = 0 \Rightarrow R''(\sigma) \neq 0.$$

Assumption \ref{Morse} implies point 1. and point 2. for $\sigma \neq \Sigma^*_{jk}$. Let us prove that $\Sigma^*_{jk}$ is a non-degenerate critical point by showing that $\Sigma^*_{jk}$ is a global maximum of $R$. Indeed, for all $\sigma$ such that $|\sigma| < 1 $: 
\begin{equation}
\left. 
\begin{aligned}  \label{ptcritique} R(\sigma) \leq R(\Sigma^*_{jk}) & \Longleftrightarrow\ \E_{\Sigma^*_{jk}}\Big(\ell(\sigma, Y)\Big) \leq \E_{\Sigma^*_{jk}}\Big(\ell(\Sigma^*_{jk}, Y)\Big) \\
& \Longleftrightarrow\ \E_{\Sigma^*_{jk}}\Big(\log \mathcal L_{jk}(\sigma, Y)\Big) \leq \E_{\Sigma^*_{jk}}\Big(\log \mathcal L_{jk}(\Sigma^*_{jk},Y)\Big)  \\
& \text{since $\ell(\sigma, y) = \log \mathcal L_{jk}(\sigma,y)$ (defined in \eqref{ell}),}\\
& \Longleftrightarrow\ \E_{\Sigma^*_{jk}}\Big(\log \dfrac{\mathcal L_{jk}(\sigma,Y)}{\mathcal L_{jk}(\Sigma^*_{jk},Y)}\Big) \leq 0.
\end{aligned}
\right\}
\end{equation}

By Jensen inequality, 
\begin{align}
\E_{\Sigma^*_{jk}}\Bigg(\log \dfrac{\mathcal L_{jk}(\sigma,Y)}{\mathcal L_{jk}(\Sigma^*_{jk},Y)}\Bigg)&  \leq \log \E_{\Sigma^*_{jk}}\Bigg(\dfrac{\mathcal L_{jk}(\sigma,Y)}{\mathcal L_{jk}(\Sigma^*_{jk},Y)}\Bigg)\nonumber\\
& = \log \int_{\R^2} \mathcal L_{jk}(\sigma,y) \mathrm{d}\mu(y) = 0, \label{Jensen}
\end{align}
since $y \mapsto \mathcal L_{jk}(\sigma, y)$ is a probability density function (see \eqref{densprob}) w.r.t. the measure $\mu$ on $\R^2$ defined in \eqref{mu}.

Equation \eqref{Jensen} implies \eqref{ptcritique}. $\Sigma^*_{jk}$ is thus a global maximum of $R$ and $R'(\Sigma^*_{jk}) = 0$. Let us prove that $R''(\Sigma^*_{jk}) \neq 0$:
\begin{align*}
 R''(\Sigma^*_{jk}) & = \E_{\Sigma^*_{jk}}\Big(\partial^2_{\sigma} \ell(\Sigma^*_{jk}, Y) \Big)\text{ by point \ref{L3.2} of Lemma \ref{L3}} \\
 & =  \E_{\Sigma^*_{jk}}\Bigg(\dfrac{\partial^2_{\sigma}\mathcal L_{jk}(\Sigma^*_{jk},Y)}{\mathcal L_{jk}(\Sigma^*_{jk},Y)} - \bigg(\dfrac{\partial_{\sigma}\mathcal L_{jk}(\Sigma^*_{jk},Y)}{\mathcal L_{jk}(\Sigma^*_{jk},Y)}\bigg)^2 \Bigg)\\
 & = - \E_{\Sigma^*_{jk}}\Bigg(\bigg(\dfrac{\partial_{\sigma}\mathcal L_{jk}(\Sigma^*_{jk},Y)}{\mathcal L_{jk}(\Sigma^*_{jk},Y)}\bigg)^2 \Bigg),
 \end{align*}
 since $\E_{\Sigma^*_{jk}}\Bigg(\dfrac{\partial^2_{\sigma}\mathcal L_{jk}(\Sigma^*_{jk},Y)}{\mathcal L_{jk}(\Sigma^*_{jk},Y)}\Bigg) = \displaystyle \int_{\R^2} \partial^2_{\sigma}\mathcal L_{jk}(\Sigma^*_{jk},y) \mathrm{d}\mu(y) = 0$ by point \ref{L3.1} of Lemma \ref{L3}. If we assume that $R''(\Sigma^*_{jk}) = 0$, then $\partial_{\sigma}\mathcal L_{jk}(\Sigma^*_{jk},Y) = 0$ a.s., which contradicts the definition of $\mathcal L_{jk}(\Sigma^*_{jk},Y)$ given in \eqref{CoupleLikelihood}.

 Accordingly, there exists $\epsilon > 0$ and $\eta > 0$ such that $R$ is $(\epsilon, \eta)$ strongly Morse. Assumption ``Morse" is satisfied.\\
 
For each couple $(j,k)$ such that $j<k$, Theorem 2 of \cite{Mei2017landscape} applied to the estimator $\widetilde \Sigma^{(n)}_{jk}$ yields: \\
Let $0 < \rho < 1$. There exists an universal constant $C_0$ such that letting $\dfrac{n}{\log n} \geq 4 C_0 \Big[ \log\Big(\dfrac{\tau(1-\delta)}{\rho} \Big) \vee 1 \Big] \Big(\dfrac{\tau^2}{\epsilon^2} \vee \dfrac{\tau^4}{\eta^2} \vee \dfrac{\tau^2L^2}{\eta^4} \Big)$ with $\tau = 2\gamma L_1$ and $L = \underset{\sigma : |\sigma| \leq 1 - \delta}{\sup} \big| R^{(3)}(\sigma) \big|$, 
$$\P\Bigg( \big|\widetilde \Sigma^{(n)}_{jk} - \Sigma^*_{jk}\big|  \leq \dfrac{2\tau}{\eta} \sqrt{C_0 \dfrac{\log n}{n} \ \Big[ \log\Big(\dfrac{\tau(1-\delta)}{\rho} \Big) \vee 1 \Big]} \Bigg) \geq 1 - \rho. $$
It follows that, for $0 < \rho < 1$ and $n$ such that $\dfrac{n}{\log n} \geq 4 C_0 \Big[ \log\Big(\dfrac{\tau(1-\delta)}{\rho} \Big) \vee 1 \Big] \Big(\dfrac{\tau^2}{\epsilon^2} \vee \dfrac{\tau^4}{\eta^2} \vee \dfrac{\tau^2L^2}{\eta^4} \Big)$, then:
\begin{align*}
\P\Bigg( \big|\big|\widetilde \Sigma^{(n)} - \Sigma^*\big|\big|_\infty  \leq \dfrac{2\tau}{\eta} \sqrt{C_0 \dfrac{\log n}{n} \ \Big[ \log\Big(\dfrac{\tau(1-\delta)}{\rho} \Big) \vee 1 \Big]} \Bigg) \geq 1 - \rho\dfrac{p(p-1)}{2},
\end{align*}
where $||A||_\infty = \underset{j,k \in \{1,\dots,p\}}{\max} |A_{jk}|$ denotes the elementwise infinite norm of the matrix $A$. This finishes the proof of Proposition \ref{Prop1}.

\subsection{Precision matrix estimator and graph recovery}

Before giving a result about the estimator of the precision matrix $\Theta^*$, let us state a third and last assumption:

\begin{enumerate}[label = \textbf{(H\arabic*)}, leftmargin = *]
\setcounter{enumi}{\value{saveenum}}
  \item  \label{Ravik} There exists some $\alpha \in ]0,1]$ such that:
 $$ \underset{e \in S^c}{\max}\ \big|\big|\Gamma_{eS}^*(\Gamma^*_{SS})^{-1}\big|\big|_1 = \big|\big|\big| \Gamma_{S^cS}^*(\Gamma^*_{SS})^{-1} \big|\big|\big|_\infty \leq 1 - \alpha,$$
 where:
 \begin{itemize}
  \item if $M \in \mathcal M_{r,m}(\R)$, $A \subset \llbracket 1,r\rrbracket$ and $B \subset \llbracket 1,m\rrbracket$, $M_{AB}$ denotes the matrix $(m_{ij})_{i \in A, j \in B}$,
  \item $S = S(\Theta^*) := E(\Theta^*) \cup \{(1,1),\dots,(p,p)\}$ where $\Theta^* = (\Sigma^*)^{-1}$ and $E(\Theta^*)$ is the set of the edges of the theoretical graph (see \eqref{edges}). Let $s := |E(\Theta^*)|$, hence $|S(\Theta^*)| = |E(\Theta^*)| + p = s + p$,
  \item $S^c = S^c(\Theta^*) = \llbracket 1,p\rrbracket^2 \setminus S(\Theta^*)$,
  \item $\Gamma^* = \Sigma^* \otimes \Sigma^*$ where $\otimes$ denotes the Kronecker matrix product. We have: $\Gamma^*_{(j,k),(l,m)} = \cov(X_j X_k, X_l X_m)$ and thus $\Gamma^*_{SS} \in \mathcal M_{s+p,s+p}(\R)$,
 \item $||u||_1 = \sum_{j=1}^d |u_j|$ for all $u \in \R^d$ is the $\ell_1$-norm,
 \item $|||U|||_\infty = \underset{j=1,\dots,d}{\max} \sum\limits_{k=1}^m |U_{jk}|$ for all $U \in \mathcal M_{d,m}(\R)$.
 \end{itemize} 
 \end{enumerate}

The underlying intuition is that this assumption \ref{Ravik} limits the influence that the
non-edge terms, indexed by $S^c$, can have on the edge-based terms, indexed by $S$ \citep{Ravikumar2011}.\\

Proposition \ref{Prop2}, set out below, provides an upper-bound for the elementwise maximum norm of the precision matrix estimator $\widetilde \Theta^{(n)}$ obtained by the graphical lasso procedure \eqref{LassoRav}. It also shows its sparsistency with respect to graphical model structure recovery. Here are some preliminary notations:
\begin{itemize}
\item $d$ is the maximum degree: 
\begin{equation} \label{degre}
d = \underset{j = 1, \dots, p}{\max} \Big|\{ k \in \llbracket 1,p \rrbracket : \Theta^*_{jk} \neq 0 \} \Big|.
\end{equation} 
\item $\kappa_{\Sigma^*}$ and $\kappa_{\Gamma^*}$ are defined by:   \begin{align} 
\label{kappaSigma}
\kappa_{\Sigma^*} & := |||\Sigma^*|||_\infty = \underset{j=1,\dots,p}{\max} \sum\limits_{k=1}^p |\Sigma^*_{jk}|,\\
\label{kappaGamma}
\kappa_{\Gamma^*} & := \bigg|\bigg|\bigg|\Big(\Gamma^*_{SS}\Big)^{-1}\bigg|\bigg|\bigg|_\infty.
\end{align}
\end{itemize}

\begin{prop} \label{Prop2} Assume \ref{Ravik} and assume that there exist some strictly positive constants $B$, $C$, $D$ and $c >2$ such that letting $\dfrac{n}{\log n} \geq C \log\big(Bp^c \big)$, we have:
\begin{align}
\P\Bigg( \big|\big|\widetilde \Sigma^{(n)} - \Sigma^*\big|\big|_\infty  \geq D \sqrt{\dfrac{\log n}{n} \log\big(Bp^c \big)}  \Bigg) \leq \dfrac{p(p-1)}{2p^c} \label{hyp},
\end{align}
where $||A||_\infty = \underset{j,k \in \{1,\dots,p\}}{\max} |A_{jk}|$ denotes the elementwise infinite norm of the matrix $A$. \\
Assume that the sample size $n$ is lower bounded as $\dfrac{n}{\log n} > D^2\log\big(Bp^c\big)\max\Big\{\dfrac{\sqrt{C}}{D}, 6(1+8\alpha^{-1})d\max\{\kappa_{\Sigma^*}\kappa_{\Gamma^*},\kappa^3_{\Sigma^*}\kappa^2_{\Gamma^*}\} \Big\}^2 $, and denote $\widetilde\Theta^{(n)}$ the unique solution of \eqref{LassoRav} and $\lambda_n = \dfrac{8D}{\alpha}\sqrt{\dfrac{\log n}{n} \log\big(Bp^c \big)}$ the regularisation parameter involved in \eqref{LassoRav}. 
Then, with probability greater than $1 - \dfrac{1}{p^{c-2}}$, we have:
\begin{enumerate}[label = (\alph*)]
\item The estimator $\widetilde \Theta^{(n)}$ of $\Theta^*$ satisfies: 
$$||\widetilde \Theta^{(n)} - \Theta^*||_\infty \leq 2D(1+8\alpha^{-1})\kappa_{\Gamma^*}\sqrt{\dfrac{\log n}{n} \log\big(Bp^c \big)}.$$
\item The estimated edges set is a subset of the true edges set: $E(\widetilde \Theta^{(n)}) \subset E(\Theta^*)$ and $E(\widetilde \Theta^{(n)})$ includes all edges $(j,k)$ with: 
$$\big|\Theta^*_{jk}\big| > 2D(1+8\alpha^{-1})\kappa_{\Gamma^*} \sqrt{\dfrac{\log n}{n} \log\big(Bp^c \big)}.$$
\end{enumerate}
\end{prop}

Proof relies on results of Theorem 1 of \cite{Ravikumar2011}, in which they study the precision matrix estimation problem in the multivariate Gaussian setting.

\begin{proof} (Proposition \ref{Prop2}) Let us check the two assumptions of Theorem 1 of \cite{Ravikumar2011}.
\vspace{0.3cm}

$\bullet$ \textbf{Incoherence assumption.} This assumption is exactly our assumption \ref{Ravik}.

\vspace{0.7cm}

$\bullet$ \textbf{Control of sampling noise.} A careful reading of \cite{Ravikumar2011} reveals that the \emph{tail conditions} of their Theorem 1 is not necessary. The required assumption, stated below, is in fact weaker, and is given in Lemma 8 of \cite{Ravikumar2011}:\\
\textit{There exists $v^* > 0$ such that for all $c > 2$ and $n$ such that $\bar\beta_f(n, p^c) \leq \dfrac{1}{v_*}$, we have:
$$ \P\Big[ ||\widetilde{\Sigma}^{(n)} - \Sigma^*||_\infty \geq  \bar\beta_f(n, p^c)\Big] \leq \dfrac{1}{p^{c-2}}, $$ 
where $\bar\beta_f(n, r) := \argmax\{\beta:f(n,\beta) \leq r\}$ for some function $f(n,\beta)$.\\
}

Setting $f(n, \beta) = \dfrac{1}{B} \exp\bigg(\dfrac{n}{\log n}\Big(\dfrac{\beta}{D}\Big)^2 \bigg)$ and $v_* = \dfrac{\sqrt{C}}{D}$ and noticing that $\dfrac{p(p-1)}{2} \leq p^2$, Assumption \eqref{hyp} gives this result. Indeed:
\begin{itemize}
\item $\bar\beta_f(n, r) = \argmax\{\beta:f(n,\beta) \leq r\} = D \sqrt{\dfrac{\log n}{n} \log\big(Br \big)}$
\item $\dfrac{n}{\log n} \geq C \log\big(Bp^c \big) \Longleftrightarrow \bar\beta_f(n, p^c) \leq \dfrac{D}{\sqrt{C}}$
\end{itemize} 
Assumption ``Control of sampling noise" is satisfied.\\

At last, let us set $\bar n_f(\beta, r) := \argmax\{n:f(n,\beta) \leq r\}$. Then, the condition $n > \bar n_f(\beta, r)$ is equivalent to $\dfrac{n}{\log n} > \log\big(Br\big)\dfrac{D^2}{\beta^2}$ since $f(n,\beta) \leq r \Longleftrightarrow \dfrac{n}{\log n} \leq \log\big(Br\big)\dfrac{D^2}{\beta^2}$.\\

We complete the proof by applying Theorem 1 of \cite{Ravikumar2011}.
\end{proof}

Finally, Propositions \ref{Prop1} and \ref{Prop2} provide the following theorem, which establishes consistency of the estimator $\widetilde \Theta^{(n)}$ in the elementwise maximum-norm:

\begin{thm} \label{Th1} Assume \ref{majoration}, \ref{Morse} and \ref{Ravik}. Let  $c > 2$, $\widetilde\Theta^{(n)}$ the unique solution of \eqref{LassoRav} and $\alpha$, $d$, $\kappa_{\Sigma^*}$ and $\kappa_{\Gamma^*}$ respectively defined in \ref{Ravik}, in \eqref{degre}, in \eqref{kappaSigma} and in \eqref{kappaGamma}. There exists some known constants $B$, $C$ and $D$ such that letting $n$ lower bounded as $\dfrac{n}{\log n} > D^2\log\big(Bp^c\big)\max\Big\{\dfrac{\sqrt{C}}{D}, 6(1+8\alpha^{-1})d\max\{\kappa_{\Sigma^*}\kappa_{\Gamma^*},\kappa^3_{\Sigma^*}\kappa^2_{\Gamma^*}\} \Big\}^2$ and $\lambda_n = \dfrac{8D}{\alpha}\sqrt{\dfrac{\log n}{n} \log\big(Bp^c \big)}$ the penalisation parameter of the equation \eqref{LassoRav}, we have, with probability greater than $1 - \dfrac{1}{p^{c-2}}$:
\begin{enumerate}[label = (\alph*)]
\item The estimator $\widetilde \Theta^{(n)}$ of $\Theta^*$ satisfies: 
$$||\widetilde \Theta^{(n)} - \Theta^*||_\infty \leq 2D(1+8\alpha^{-1})\kappa_{\Gamma^*}\sqrt{\dfrac{\log n}{n} \log\big(Bp^c \big)}.$$
\item $E(\widetilde \Theta^{(n)}) \subset E(\Theta^*)$ and $E(\widetilde \Theta^{(n)})$ includes all edges $(j,k)$ with: 
$$\big|\Theta^*_{jk}\big| > 2D(1+8\alpha^{-1})\kappa_{\Gamma^*} \sqrt{\dfrac{\log n}{n} \log\big(Bp^c \big)}.$$
In other words, the graph structure of latent Gaussian encoded in $\Theta^*$ is consistently recovered as long as: $\big|\Theta^*_{jk}\big| > 2D(1+8\alpha^{-1})\kappa_{\Gamma^*} \sqrt{\dfrac{\log n}{n} \log\big(Bp^c \big)}.$
\end{enumerate}
\end{thm}

The parameter $c$ of Theorem \ref{Th1} is a user-defined parameter. The larger $c$ is, the larger the probability for which results of Theorem \ref{Th1} hold is. However, large values of this parameter lead to more stringent requirements on the sample size $n$.

\section{Simulation studies}

\subsection{Simulation settings}

In almost all of this Section (unless otherwise stated), we use the following simulation settings. We simulate $n = 500$ observations of a $p = 100$-Gaussian vector $X$ centered and scaled. Graph structure is a chain, that is $X_1 \longleftrightarrow X_2 \longleftrightarrow \dots \longleftrightarrow X_{100}$. The graph contains then 99 edges. Data have been simulated with the R function \texttt{huge.generator}, option \texttt{graph = "band"} of the package \texttt{huge}.\\

Two different settings of the points of truncation are presented:
\begin{itemize}[label = $\bullet$]
\item identical points of truncation: $a = -0.5$ and $b = 2$,
\item decreasing points of truncation: $a = -1$, $b = \texttt{seq(2,0.5, length = p)}$.\\
\end{itemize}

We then apply the estimation procedure described in Subsection \ref{ssect:EstiMat}. We assume that the points of truncation are known and we use the estimators $\widetilde \Sigma^{(n)}$ and $\widetilde \Theta^{(n)}$ of the covariance and precision matrices, respectively defined at Step 1 \eqref{EstiSigma} and at Step 2 \eqref{LassoRav}. \\

For the estimation of the precision matrix with the graphical lasso procedure, we use the function \texttt{huge} of the same package, option \texttt{method = "glasso"}. Unfortunately, theoretical results of Section \ref{sect:ResultatsTheo} do not give an explicit choice of the penalty parameter. We thus choose the penalty parameter with ``stars" and ``ebic" methods, implemented in the package \texttt{huge}.\\

Simulations address several problems. Let first explicit the two different procedures used for these simulations, respectively called ``our procedure" and ``graphical lasso directly on truncated data" (shortened in ``Glasso"). The first one is our procedure, which consists in replacing the empirical covariance matrix of $X$ in the graphical lasso by our estimator $\widetilde \Sigma^{(n)}$. The second one is the graphical lasso directly applied to the truncated data, which consists in replacing the empirical covariance matrix of $X$ in the graphical lasso by the empirical covariance matrix of the truncated vector $Y$. Here are the problems addressed in the following subsections: 
\begin{itemize}
\item Efficiency of our procedure. Does using our estimator for the covariance matrix really improve graph estimation? 
\item Impact of the points of truncation, that is how the values of the truncation points impact edges detection.
\item Is our procedure as efficient with other graph structures?
\end{itemize}

To study efficiency of these procedures on graph estimation, we make 50 i.i.d. repetitions of the procedure and we represent the detection rates of each of the $\dbinom{100}{2} = 4950$ potential edges.

\subsection{Efficiency}\label{ssect:efficacite}

In this subsection, we aim at illustrating the efficiency of our procedure. For that, we compare detection rates of each potential edges obtained with our procedure and with graphical lasso directly on truncated data. \\

\begin{figure}%
\centering
\subfigure[Our procedure, identical points of truncation: $a = -0.5$ and $b = 2$.]{\includegraphics[scale = 0.42]{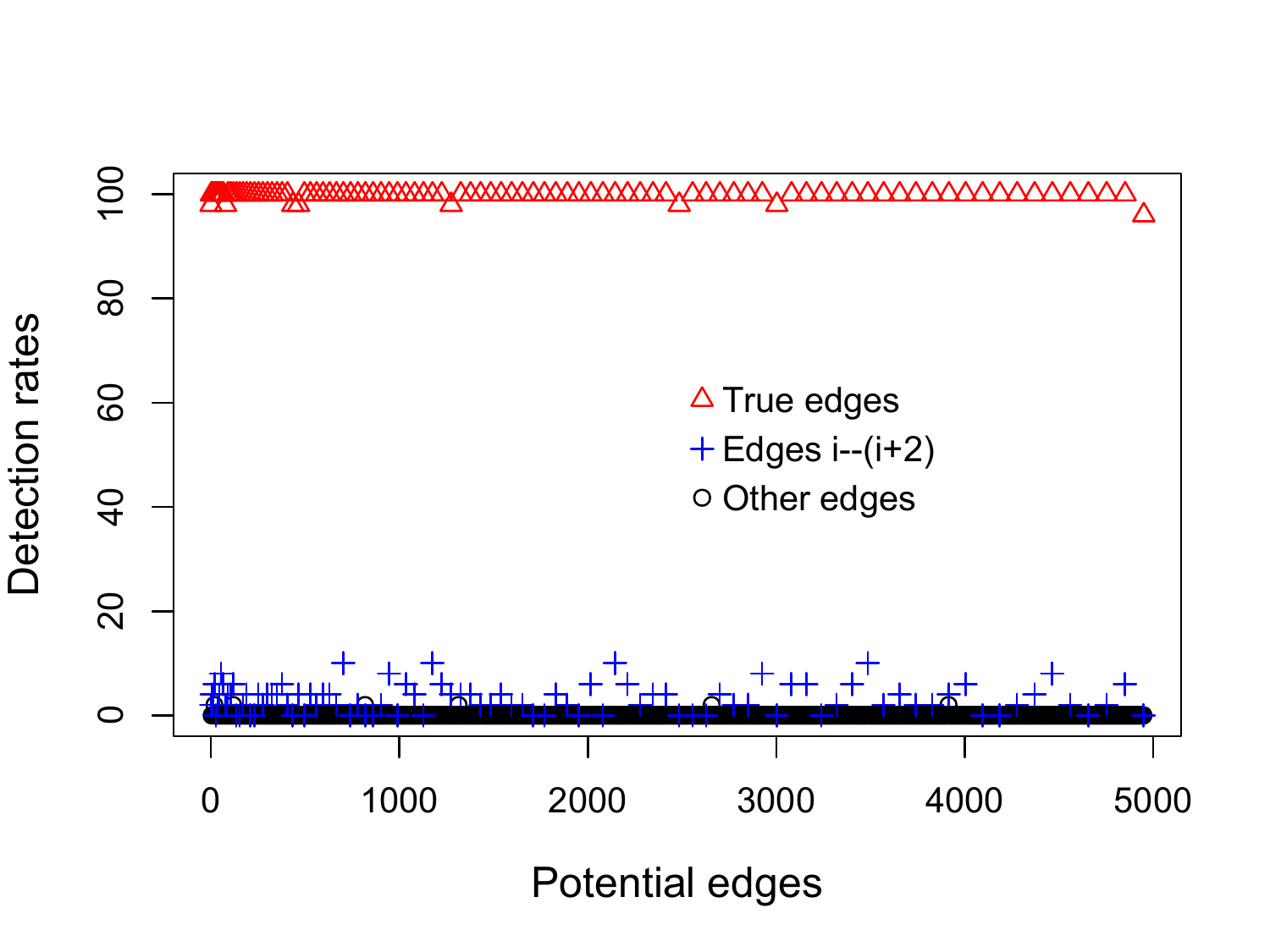}\label{fig:eff1_us}}\qquad
\subfigure[Glasso, identical points of truncation: $a = -0.5$ and $b = 2$.]{\includegraphics[scale = 0.42]{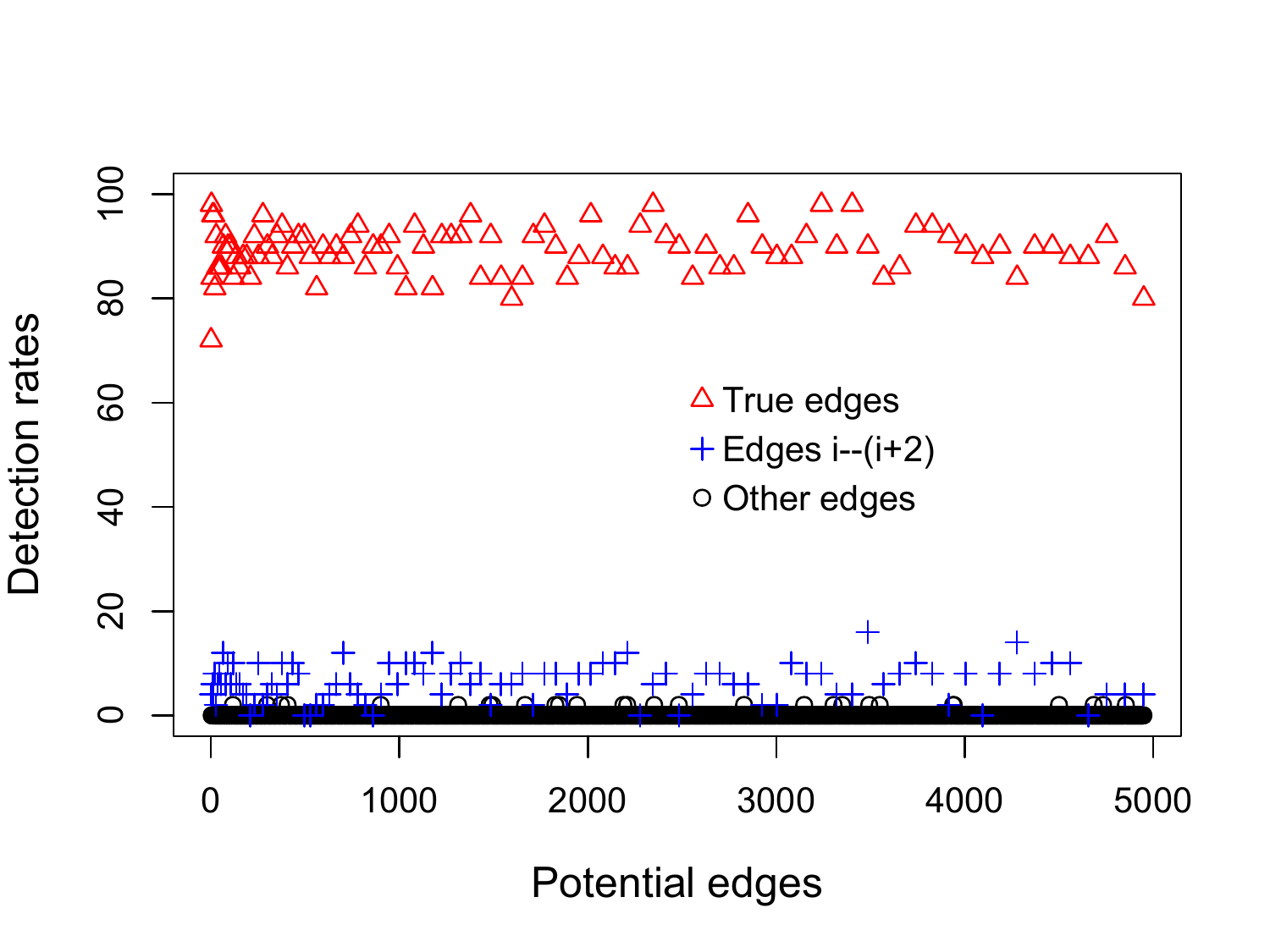}}\\
\caption{Comparison of detection rates obtained by our procedure and by graphical lasso directly on truncated data. Identical points of truncation setting. Detection rates are obtained on 50 i.i.d. repetitions for $n = 500$ observations of $p = 100$ variables. True edges are represented with red triangles, (false) edges of type $X_i - X_{i+2}$ with blue crosses and other false edges with black circles.}\label{fig:eff1}
\end{figure}

\begin{figure}%
\centering
\subfigure[Our procedure, decreasing points of truncation: $a = -1$, $b = \texttt{seq(2,0.5, length = p)}$.]{\includegraphics[scale = 0.42]{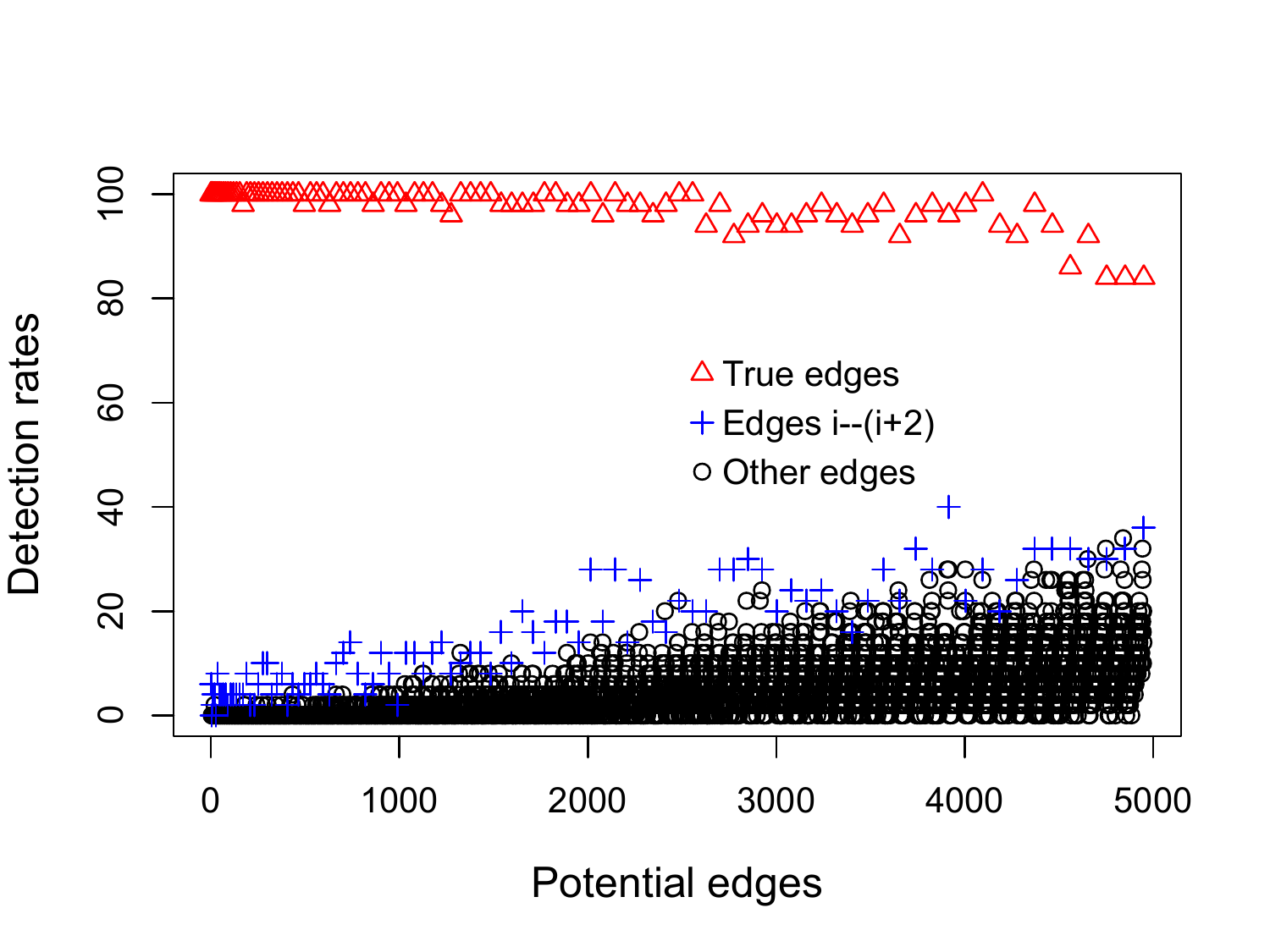}\label{fig:eff2_us}}\qquad
\subfigure[Glasso, decreasing points of truncation: $a = -1$, $b = \texttt{seq(2,0.5, length = p)}$.]{\includegraphics[scale = 0.42]{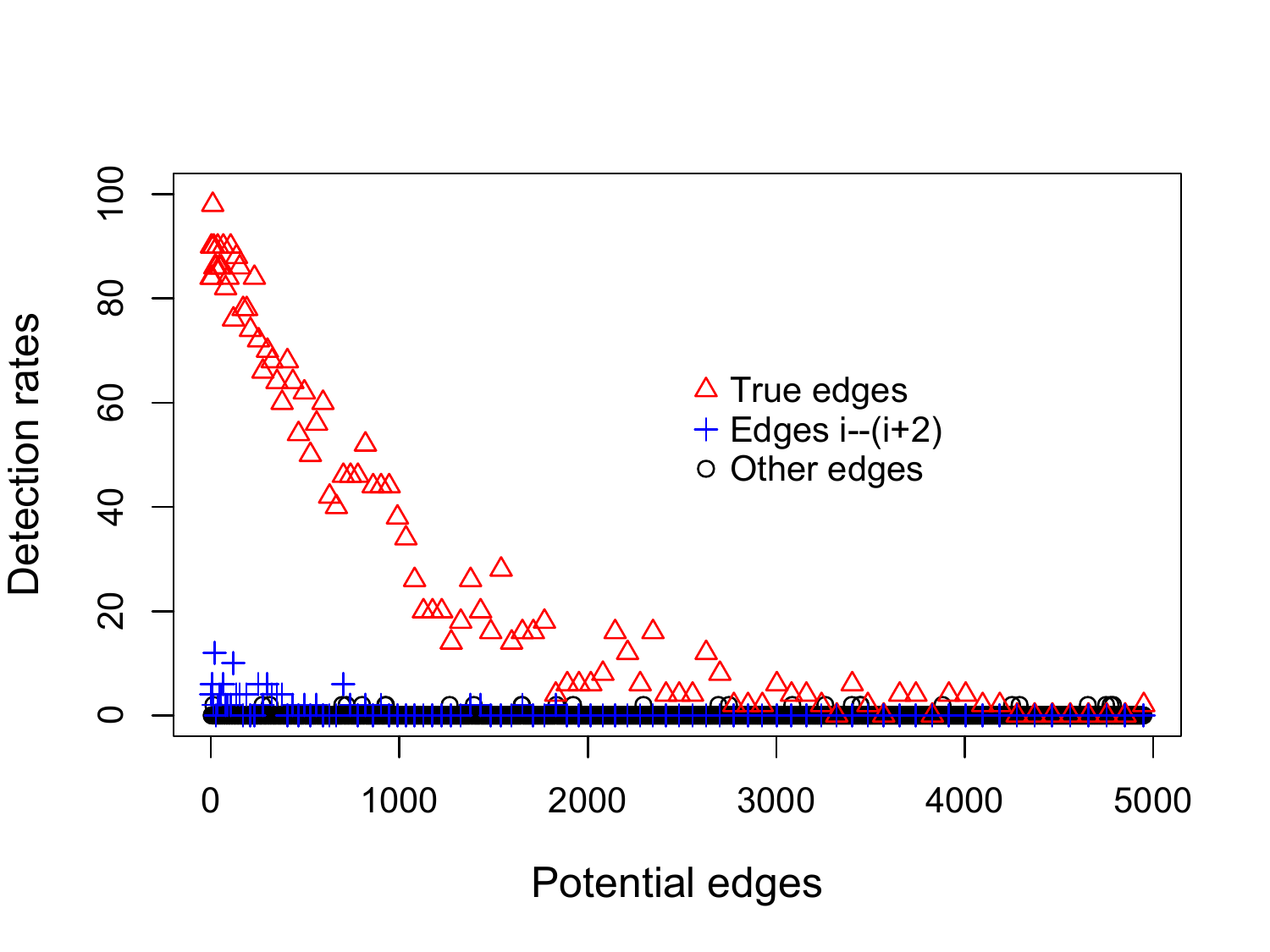}}\\
\caption{Comparison of detection rates obtained by our procedure and by graphical lasso directly on truncated data. Decreasing points of truncation setting. Detection rates are obtained on 50 i.i.d. repetitions for $n = 500$ observations of $p = 100$ variables. True edges are represented with red triangles, (false) edges of type $X_i - X_{i+2}$ with blue crosses and other false edges with black circles.}\label{fig:eff2}
\end{figure}

\noindent
\textbf{Results and comments.} Figures \ref{fig:eff1} and \ref{fig:eff2} illustrate detection rates for our procedure and for graphical lasso directly on truncated data $Y$. We represent detection rates of each potential edge by unfolding the matrix with \texttt{upper.tri}, which implies that the first potential edge whose detection rate is displayed is $X_2 \longleftrightarrow X_1$, then $X_3 \longleftrightarrow X_1$, $X_3 \longleftrightarrow X_2$, $X_4 \longleftrightarrow X_1$, $\dots$, $X_4 \longleftrightarrow X_3$ and so on. The 99 true theoretical edges, that is edges $X_i \longleftrightarrow X_{i+1}$, are displayed with red triangles and in the following order $X_1 \longleftrightarrow X_2$, $X_2 \longleftrightarrow X_3$, ..., $X_{99} \longleftrightarrow X_{100}$. Detection rates of edges of type $X_i \longleftrightarrow X_{i+2}$, which are not true edges, are displayed with blue crosses. We distinguish these edges because these interactions can be relatively strong because of the indirect link through $X_{i+1}$.  \\
Figure \ref{fig:eff1} shows identical points of truncation setting. We can observe that our method gives better results: true edges are better detected and other edges are less detected than with Glasso. For example, true edges are all detected more than 96\% whereas Glasso detects them between 80\% and 100\%: 66 of them are detected at most 90\%. Edges of type $X_i \longleftrightarrow X_{i+2}$ tend to be more detected with Glasso.\\
Figure \ref{fig:eff2} shows decreasing points of truncation setting. Efficiency of our procedure is even more convincing in this setting. Indeed, true edges are much better detected by our procedure (more thant 80\% whereas only 60 true edges are detected at most 80\% by Glasso). These differences are probably due to the zero rate in the truncated data $Y$ which grows from 20\% to 50\% according to the variables. Consequently, edges between variables whose zeros rate is high (that is edges $X_i \longleftrightarrow X_{i+1}$ for $i$ close to 100) tend to be less detected. The phenomenon is even stronger with Glasso. Besides, the other edges (the false ones) tend to be slightly less detected with Glasso. For our procedure, the false edges are more detected when the truncation points induce a high rate of zero for the involved variables.\\

\subsection{Impact of points of truncation}

To expand this section, we illustrate the impact of the points of truncation values. For that, we briefly compare results obtained for both ``identical" and ``decreasing" settings. Observations of the underlying Gaussian vector are the same and we only change values of the points of truncation according to the chosen setting.\\

\begin{figure}%
\centering
\subfigure[Points of truncation setting: $a = -1$ and $b = 1$.]{\includegraphics[scale = 0.42]{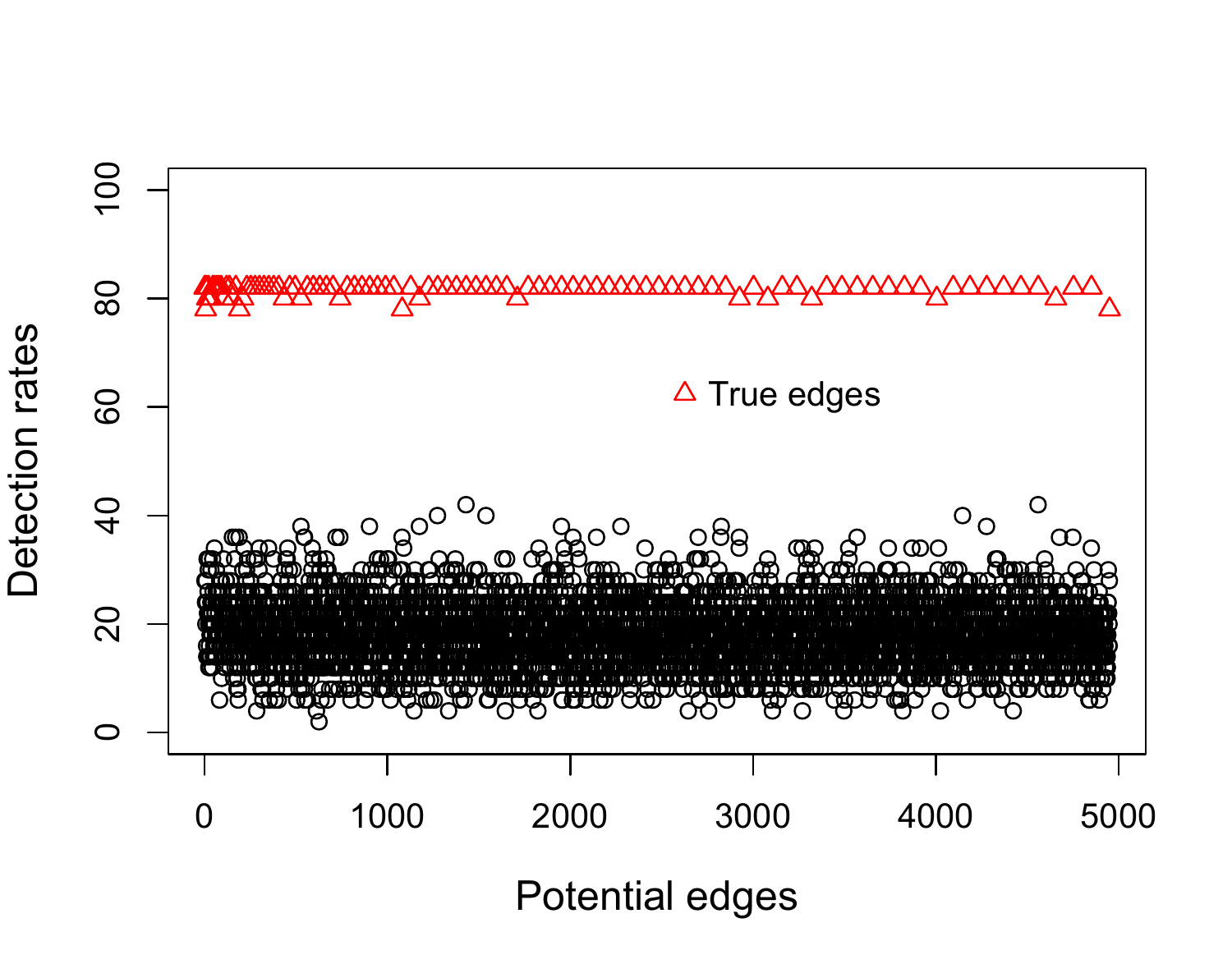}}\qquad
\subfigure[Identical points of truncation: $a = -0.5$ and $b = 2$.]{\includegraphics[scale = 0.42]{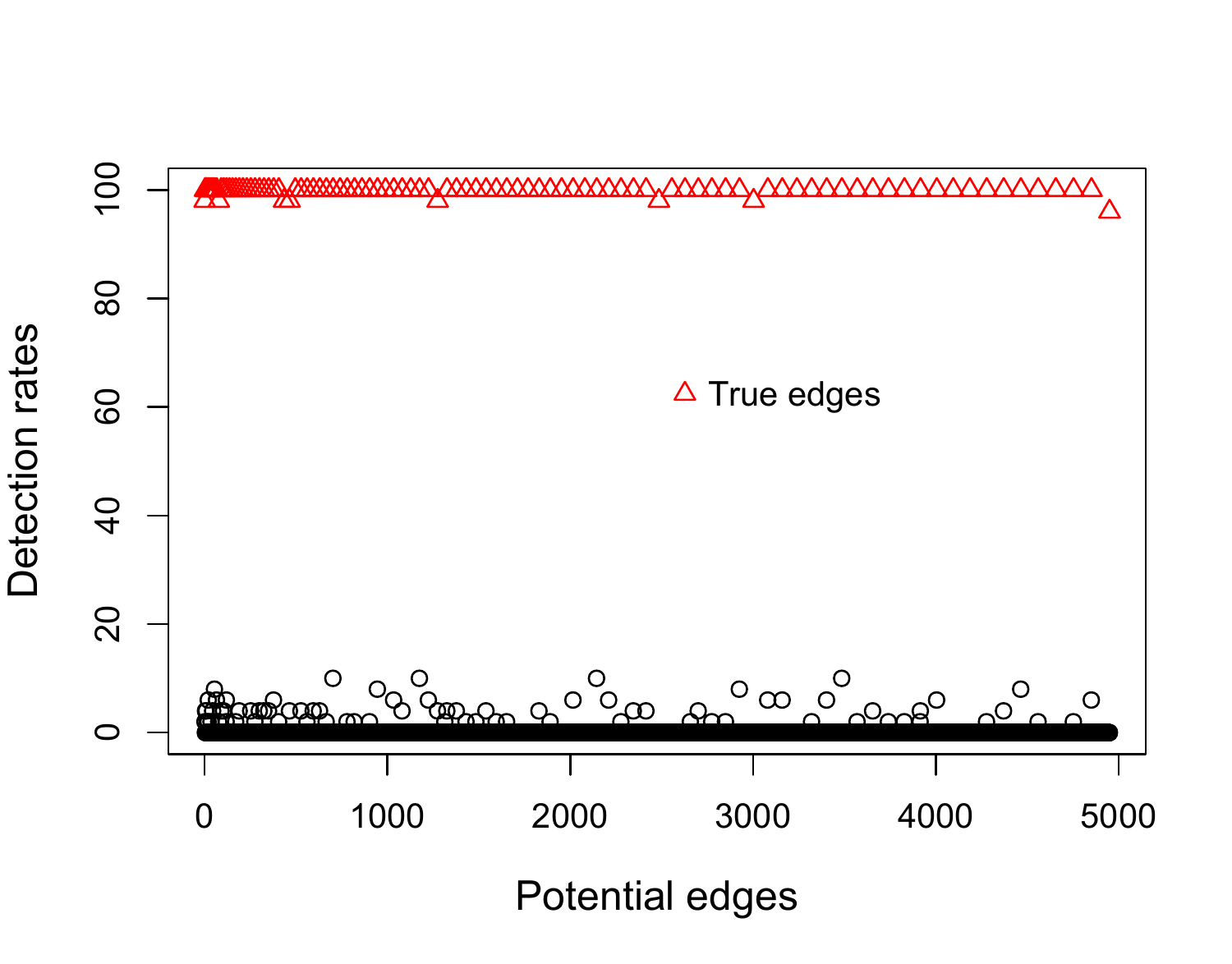}}\\
\caption{Comparison of detection rates for two points of truncation settings for which the zero inflation is similar (around 32\%). Detection rates are obtained on 50 i.i.d. repetitions for $n = 500$ observations of $p = 100$ variables. True edges are represented with red triangles.}\label{fig:Queue}
\end{figure}

\noindent
\textbf{Results and comments.} Figures \ref{fig:eff1_us} and \ref{fig:eff2_us} respectively show detection rates for ``identical" and ``decreasing" points of truncation settings. \\
The zero inflation of truncated data $Y$ of the first setting is around 33\%. Zero inflation of truncated data of the ``decreasing" setting decreases from 20\% (for $Y_1$) to 50\% (for $Y_{100}$). Detection rates of potential edges are represented by unfolding the matrix with \texttt{upper.tri}, which implies that the first potential edge whose detection rate is displayed is $X_2 \longleftrightarrow X_1$, then $X_3 \longleftrightarrow X_1$, $X_3 \longleftrightarrow X_2$, $X_4 \longleftrightarrow X_1$, $\dots$, $X_4 \longleftrightarrow X_3$ etc. The 99 true theoretical edges, that is the edges $X_i \longleftrightarrow X_{i+1}$, are displayed with red triangles and in the following order $X_1 \longleftrightarrow X_2$, $X_2 \longleftrightarrow X_3$, $\dots$, $X_{99} \longleftrightarrow X_{100}$. Thus, we can observe that edges involving variables whose zero inflation is close to 50\% have worse detection rates: true edges are less detected whereas false edges have higher detection rates.\\
In short and as expected, zero inflation impacts detection rates: the more the zero inflation is, worse is the detection rate.\\
An other phenomenon, not noticeable in Figures \ref{fig:eff1_us} and \ref{fig:eff2_us} occurs. This phenomenon is also linked to zero inflation and is noticeable in Figure \ref{fig:Queue}. It points out that detection rates does not only depend on zero inflation but also on the observation window. Indeed, Figure \ref{fig:Queue} exhibits detection rates for ``identical" points of truncation setting, that is $a = -0.5$ and $b = 2$, and for an other setting  $a = -1$ and $b = 1$. For both settings, zero inflation is about 32\%. However, results obtained for ``identical" points of truncation are much better. The underlying idea is that the observation of the Gaussian variable gives more informations between $-0.5$ and $2$ than between $-1$ and $1$, especially for covariances estimation.

\subsection{Other graph structures}

\begin{figure}%
\centering
\subfigure[``Random" graph structure.]{\includegraphics[scale = 0.42]{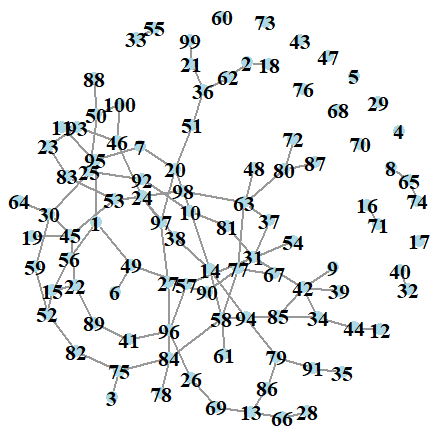}}\qquad
\subfigure[``Hub" graph structure.]{\includegraphics[scale = 0.42]{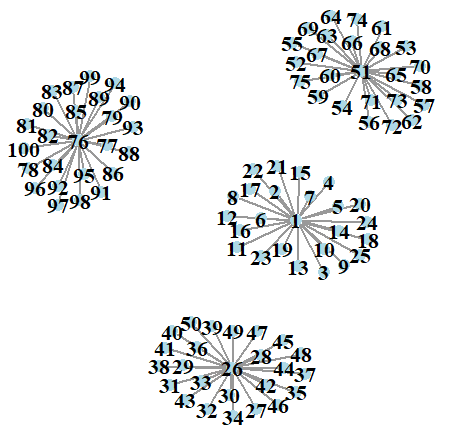}}\\
\caption{Graphical representation of the two graphs used in this subsection: ``random" and ``hub".}\label{fig:graphs_shape}
\end{figure}

%
Previously, we restrict to only one graph structure, the chain structure which tends to give satisfactory results in general. To fulfil these simulation studies, we present some results with other graph structures:
\begin{itemize}
\item The ``random" structure. There exists an edge between two variables with probability 1/50. Data have been simulated with the R function \texttt{huge.generator}, options \texttt{graph = 'random', prob = 1/50}. Resulting graph has $103$ edges.
\item The ``hub" structure. Variables are split into 4 groups of 25. Inside each group, one of the variable is a ``hub" and is connected to all the variables of its group. Data have been simulated with the R function \texttt{huge.generator}, options \texttt{graph = 'hub', g = 4}. Resulting graph has $96$ edges.
\end{itemize}
These graphs are displayed in Figure \ref{fig:graphs_shape}. Points of truncation are set to identical ($a = -0.5$ and $b = 2$). We simulate $n = 500$ observations of $p = 100$ variables and we compare detection rates of each potential edge obtained with our procedure and with graphical lasso directly on truncated data like in Subsection \ref{ssect:efficacite}. \\

\begin{figure}%
\centering
\subfigure[``Random" graph structure.]{\includegraphics[scale = 0.45]{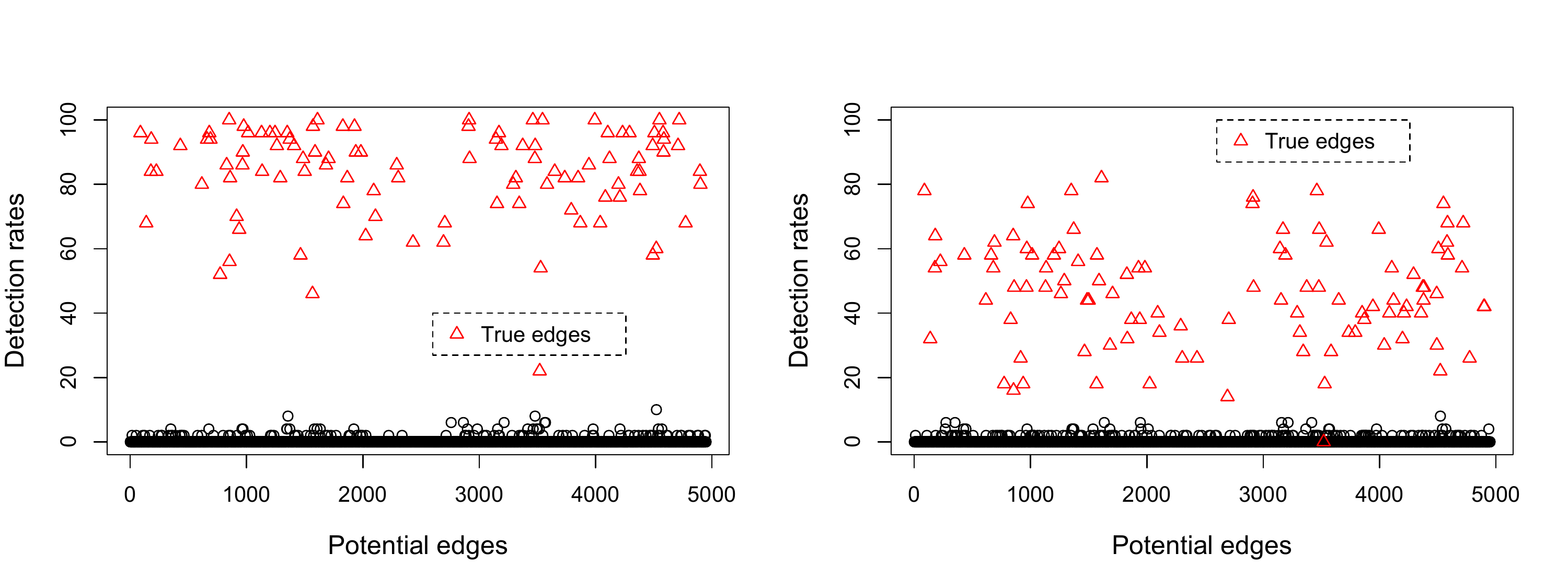} \label{fig:random}}\qquad
\subfigure[``Hub" graph structure.]{\includegraphics[scale = 0.45]{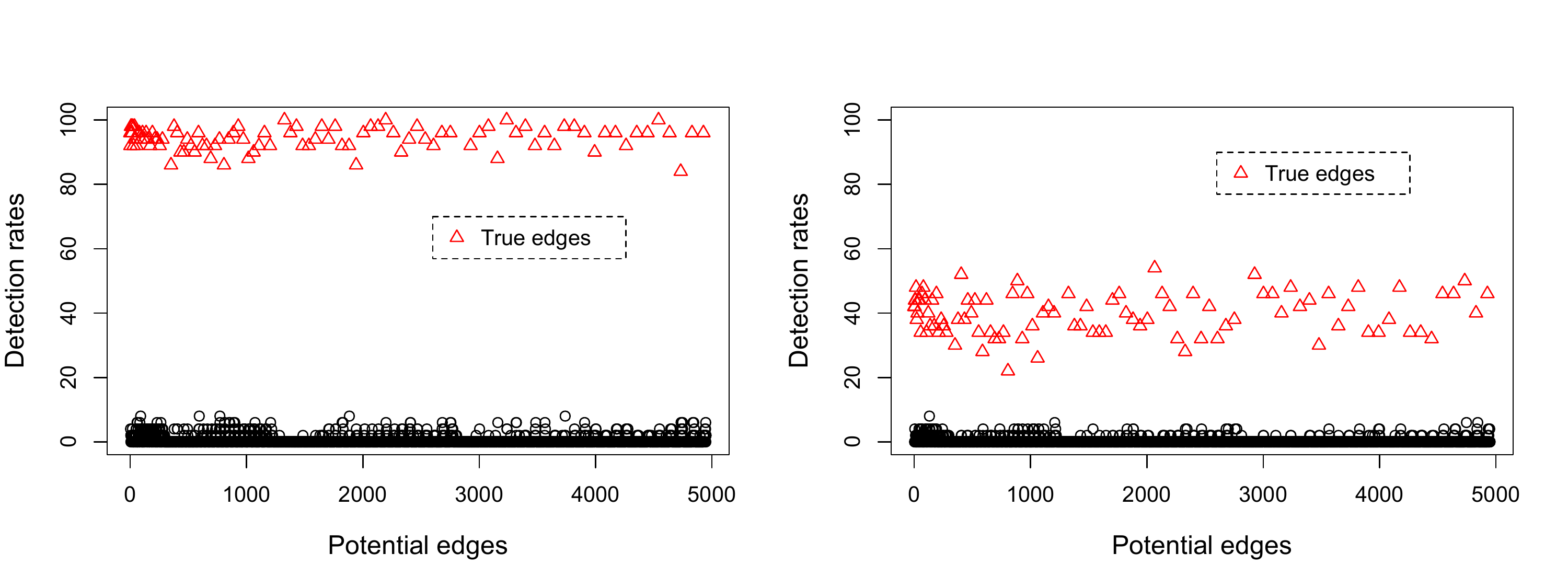}}\\
\caption{Comparison of detection rates obtained by our procedure (left) and by graphical lasso (right) directly on truncated data. Identical points of truncation setting. ``random" and ``hub" graph structures. Detection rates are obtained on 50 i.i.d. repetitions for $n = 500$ observations of $p = 100$ variables. True edges are represented with red triangles.}\label{fig:structures}
\end{figure}

\begin{figure}%
\centering
\subfigure[``Random" graph structure.]{\includegraphics[scale = 0.42]{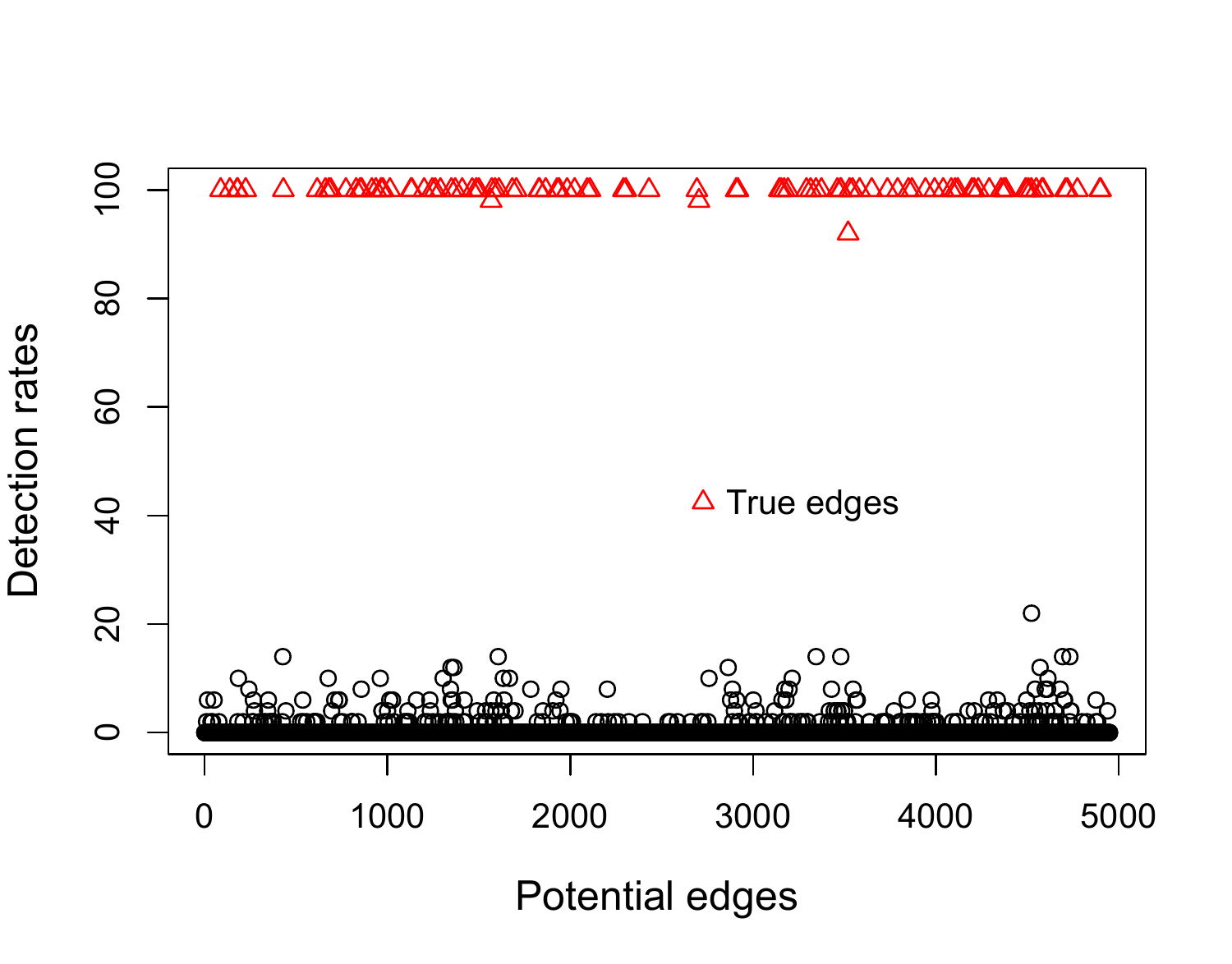}\label{fig:notrunc_random}}\qquad
\subfigure[``Hub" graph structure.]{\includegraphics[scale = 0.42]{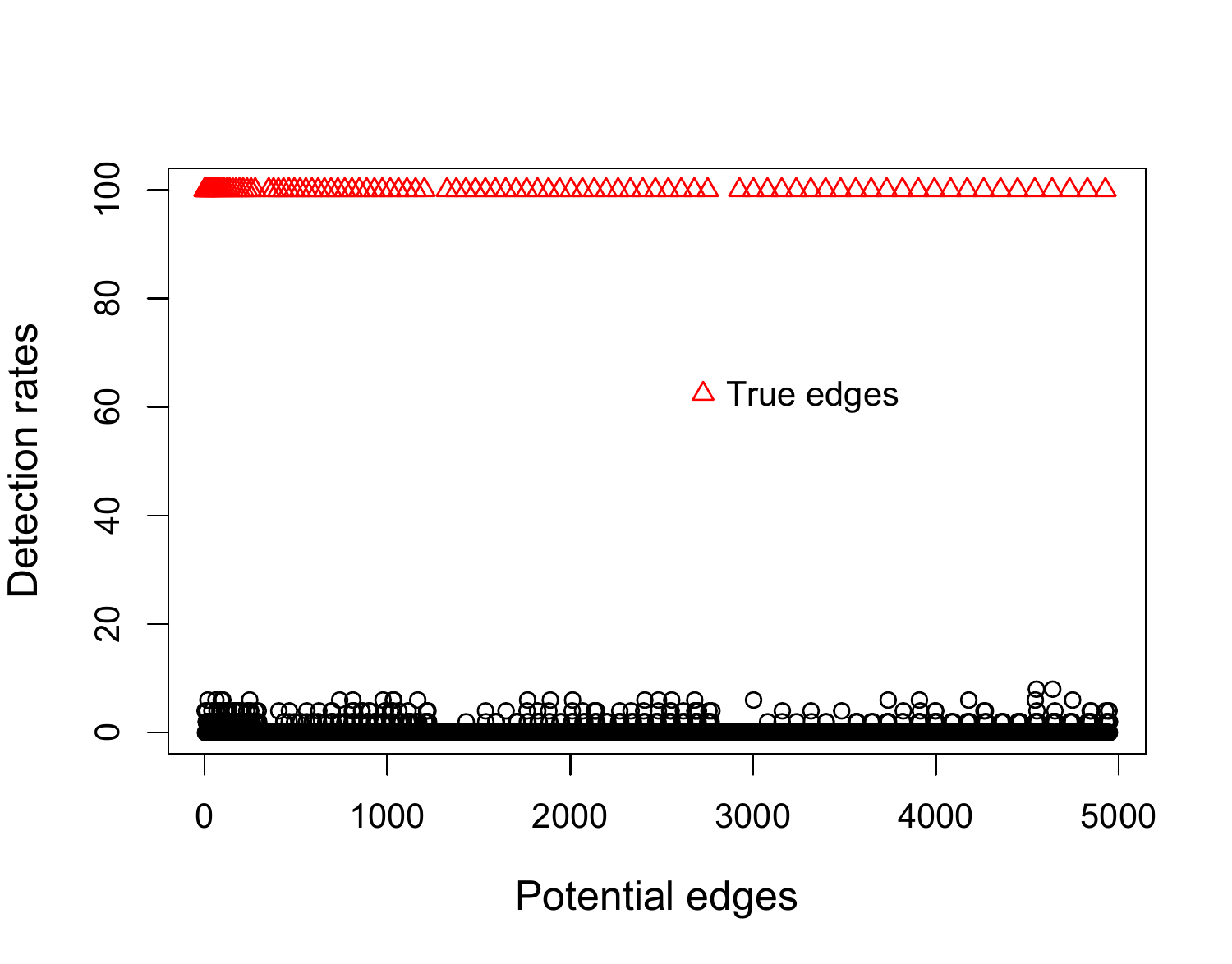}\label{fig:notrunc_hub}}\\
\caption{Detection rates obtained with graphical lasso on Gaussian data $X$. ``random" and ``hub" graph structures. Detection rates are obtained on 50 i.i.d. repetitions for $n = 500$ observations of $p = 100$ variables. True edges are represented with red triangles.}\label{fig:notrunc}
\end{figure}

\noindent
\textbf{Results and comments.} Figure \ref{fig:structures} displays these comparisons for ``random" and ``hub" graph structures. \\
Results obtained with our procedure are indeed slightly less satisfying than with ``chain" structure. Detection rates with ``hub" graph structure are a little better than with ``random". For ``random", all the false edges are detected less than 10\% and true edges are all detected more than 46\% except for the edge $X_{34} \longleftrightarrow X_{85}$ (22\%). This edge is also slightly less detected (92\% whereas at least 98\% for the other true edges) when we directly apply graphical lasso on untruncated Gaussian data $X$ (see Figure \ref{fig:notrunc_random}). Thus, this low detection rate is not only due to our procedure but probably also to the data itself. ``Hub" structure gives good results: false edges are detected at most 8\% and true edges at most 84\%, which is also the case when we apply graphical lasso on untruncated data $X$ (see Figure \ref{fig:notrunc_hub}).\\
In comparison, Glasso (on truncated data $Y$) always gives less good results. For ``random", the edge $X_{34} \longleftrightarrow X_{85}$ is for example never detected and the other true edges are detected between 14\% and 82\%. ``Hub" structure exhibits the most striking difference: true edges are detected at most 54\% (but at least 22\%).

\section{Discussion}

In this paper, we proposed a procedure for graph estimation in a zero-inflated Gaussian model. In this model, zero-inflation is obtained by double truncation (left and right) of Gaussian data. More precisely, the goal is to retrieve the underlying graph structure given by the precision matrix of the Gaussian data with the doubly truncated data. Our procedure includes two steps: the first one consists in estimating the covariance matrix terms to terms by maximising the corresponding bivariate marginal log-likelihood of the truncated vector. The second one relies on the graphical lasso procedure to obtain a lasso estimation of the precision matrix. We then proved some theoretical convergence guarantees with regard to graph estimation. The first result states rate convergence about the covariance matrix estimator. The second one provides sparsistency of the precision matrix estimator with respect to graph structure recovery.\\
Practically, simulations studies also corroborate efficiency of our procedure. They also show that our procedure is more appropriate than using graphical lasso directly on truncated data (without a preliminary estimation of the covariance matrix).\\ 
However, results depend on graph structure and simulations exhibit well that some graph structures are more favorable to graph recovery. \\

This work only deals with a double truncation. Yet, our procedure seems to be practically efficient with a single truncation (right or left) but proof of our theoretical results requires the both right and left points of truncation and they do not hold in the unilateral setting. It would be interessant to address this case later, perhaps using different tools.

\appendix

\section{Proof of Lemma \ref{L1}}

\begin{proof} (Lemma \ref{L1}) It is sufficient to show the existence of such a constant $\gamma_{jk} > 0$ for $j < k$ fixed
. Let $j < k$, $\sigma \in [-1 + \delta, 1 - \delta]$ and $(y_j, y_k) \in [a_j,b_j]\times[a_k,b_k]$. The proof naturally falls into four parts:

$\bullet$ $a = b = 1$:  \begin{align} \label{phi11}
\phi_{11,jk}(\sigma, y_j, y_k) 
& = \dfrac{1}{2\pi\sqrt{1 - \sigma^2}} \exp\bigg[ - \dfrac{y_j^2 - 2\sigma y_jy_k + y_k^2}{2(1-\sigma^2)} \bigg].
\end{align} 
Since $(y_j,y_k) \in [a_j,b_j] \times [a_k,b_k]$ and $\delta^2 \leq 1- \sigma^2 \leq 1$, $\phi_{11,jk}$ is continuous and strictly positive on a compact of $\R^3$.

\vspace{0.2cm}
$\bullet$ $a = 0, b = 1$: An easy computation yields that
\begin{align}
\phi_{01,jk}(\sigma, y_k) 
& = \dfrac{1}{\sqrt{2\pi}}\exp\bigg(-\dfrac{y_k^2}{2}\bigg)\bigg[1 - F\Big(\dfrac{b_j-\sigma y_k}{\sqrt{1-\sigma^2}}\Big) + F\Big(\dfrac{a_j-\sigma y_k}{\sqrt{1-\sigma^2}} \Big) \bigg], \label{phi01}
\end{align}
where $F$ denotes the c.d.f. of $\mathcal N(0,1)$. Since $y_k \in [a_k,b_k]$, $-\infty < a_k < b_k < \infty$ and $\delta^2 \leq 1- \sigma^2 \leq 1$, $\phi_{01,jk}$ is hence continuous and strictly positive on a compact of $\R^2$.

\vspace{0.2cm}
$\bullet$ $a = 1, b = 0$: Analogous to $a = 0, b = 1$.

\vspace{0.2cm}
$\bullet$ $a = b = 0$: 
\begin{align*}
\phi_{00,jk}(\sigma) 
& = \int_{[a_k,b_k]^c} \dfrac{1}{\sqrt{2\pi}}\exp\bigg(-\dfrac{y^2}{2}\bigg)\bigg[1 - F\Big(\dfrac{b_j-\sigma y}{\sqrt{1-\sigma^2}}\Big) + F\Big(\dfrac{a_j-\sigma y}{\sqrt{1-\sigma^2}} \Big) \bigg] \mathrm{d}y.
\end{align*}
Since $\delta^2 \leq 1- \sigma^2 \leq 1$, $\phi_{00,jk}$ is hence continuous and strictly positive on a compact of $\R$.
\end{proof}

\section{Proof of Lemma \ref{L2}}

\begin{proof} (Lemma \ref{L2}) In the same manner as for the proof of Lemma \ref{L1}, it is sufficient to show the result for $j < k$ fixed.  Fix $j < k$ and show that for all $a, b \in \{0,1\}$, the function $\phi_{ab,jk}$ is $C^3$ on the compact $[-1+\delta,1-\delta]\times[a_j,b_j]\times[a_k,b_k]$ (it is even $C^\infty$). Hence, for all $m \in \{1,2,3\}$, $\partial^m_\sigma \phi_{ab,jk}$ is continuous on a compact of $\R^3$, which establishes the result. We naturally distinguish the four cases (see \eqref{phi01} for expressions of $\phi_{ab,jk}$):

$\bullet$ $a = b = 1$: 
$\phi_{11,jk}$ is $C^3$ on $]-1,1[ \times \R^2$ hence on $[-1 + \delta, 1-\delta] \times [a_j,b_j] \times [a_k,b_k]$.

\vspace{0.2cm}
$\bullet$ $a = 0, b = 1$:  
Since $F$ is $C^\infty$ hence $C^3$ on $\R$, $\phi_{01,jk}$ is $C^3$ on 
$[-1 + \delta, 1-\delta] \times [a_k,b_k]$.

\vspace{0.2cm}
$\bullet$ $a = 1, b = 0$: Analogous to $a = 0, b = 1$.

\vspace{0.2cm}
$\bullet$ $a = b = 0$ : 
Let $m \in \{1,2,3\}$. We apply Lebesgue theorem for continuity and differentiability of integrals with parameters: 
\begin{itemize}[leftmargin = 2 cm]
\item $\sigma \mapsto \phi_{01,jk}(\sigma, y)$ is $C^3$ on $[-1+\delta,1-\delta]$.
\item Straigthforward calculations of derivatives of $\phi_{01,jk}(\sigma, y)$ w.r.t. $\sigma$ show that, for $\sigma \in [-1+\delta, 1-\delta]$:
$$ \Big|\partial^m_\sigma \phi_{01,jk}(\sigma, y)\Big| \leq C(a_j, b_j, a_k, b_k, \delta, m)\exp\bigg(- \dfrac{y^2}{2}\bigg),$$ where $C(a_j, b_j, a_k, b_k, \delta, m)$ is a positive constant depending on $a_j, b_j, a_k, b_k, \delta$ and $m$ and $y  \mapsto C(a_j, b_j, a_k, b_k, \delta, m)\exp\bigg(- \dfrac{y^2}{2}\bigg)$ is integrable on $[a_k,b_k]^c$.
\end{itemize}
It follows that $\phi_{00,jk}$ is $C^3$ on $[-1 + \delta, 1-\delta]$.
\end{proof}

\section{Proof of Lemma \ref{L3}}

\begin{proof} (Lemma \ref{L3}) 
\begin{enumerate}[leftmargin = *]
\item Let $l \in \mathbb N^*$. First, for all $\sigma \in [-1+\delta, 1-\delta]$, we have:

\begin{equation}\label{densprob}
\begin{aligned}
\displaystyle \int_{\R^2} \mathcal L_{jk}(\sigma, y) \mathrm{d}\mu(y) 
= & \iint_{\R^2} f(x,y,\sigma) \mathrm{d}x\mathrm{d}y = 1.
\end{aligned}
\end{equation}

It remains to prove that:
\begin{align} \label{fri}
\partial^l_{\sigma}\Bigg(\displaystyle \int_{\R^2} \mathcal L_{jk}(\sigma, y) \mathrm{d}\mu(y)\Bigg) = &\ \partial^l_{\sigma}\phi_{00,jk}(\sigma) + \int_{a_k}^{b_k} \partial^l_{\sigma}\phi_{01,jk}(\sigma, y) \mathrm{d}y + \int_{a_j}^{b_j} \partial^l_{\sigma}\phi_{10,jk}(\sigma, x) \mathrm{d}x \hspace{2cm}\\
& + \iint_{[a_j,b_j]\times[a_k,b_k]} \partial^l_{\sigma}\phi_{11,jk}(\sigma, x,y) \mathrm{d}x\mathrm{d}y.\nonumber
\end{align}

Let us deal with each of these terms:

\vspace{0.2cm}
$\bullet$ For $a = 0, b = 0$: it is obvious. 

For the following terms, we use Lebesgue theorem for differentiability of integrals with parameters.

\vspace{0.2cm}
$\bullet$ For $a = 0, b = 1$ (and $a = 1, b = 0$): According to \eqref{phi01}, it is clear that $\phi_{01,jk}$ is $C^\infty$ on the compact $[1-\delta, 1+\delta]\times[a_k,b_k]$, which establishes the formula.

\vspace{0.2cm}
$\bullet$ For $a = 1, b = 1$: Analogously, \eqref{phi11} shows that $\phi_{11,jk}$ is $C^\infty$ on the compact $[1-\delta, 1+\delta]\times[a_j,b_j]\times[a_k,b_k]$.

\item Let $l \in \mathbb N^*$. Let us first clarify some notations:
\begin{align*}
\E_{\Sigma^*_{jk}}\Big(\log \mathcal L_{jk}(\sigma, Y)\Big) = &\ \E_{\Sigma^*_{jk}}\Big(\ell(\sigma, Y)\Big) = R(\sigma)\\
= & \int_{\R^2} \log \mathcal L_{jk}(\sigma, y) \mathcal L_{jk}(\Sigma^*_{jk}, y)\mathrm{d}\mu(y) \\
= &\ \phi_{00,jk}(\Sigma^*_{jk})\log\phi_{00,jk}(\sigma) + \int_{a_k}^{b_k} \phi_{01,jk}(\Sigma^*_{jk}, y) \log \phi_{01,jk}(\sigma, y) \mathrm{d}y \\
& + \int_{a_j}^{b_j} \phi_{10,jk}(\Sigma^*_{jk}, x) \log \phi_{10,jk}(\sigma, x) \mathrm{d}x \\
& + \iint_{[a_j,b_j]\times[a_k,b_k]} \phi_{11,jk}(\Sigma^*_{jk}, x,y) \log \phi_{11,jk}(\sigma, x,y) \mathrm{d}x\mathrm{d}y.
\end{align*}
On the other hand, 
\begin{align*}
\hspace{-0.5cm}\E_{\Sigma^*_{jk}}\Big(\partial_{\sigma}^l \log \mathcal L_{jk}(\sigma, Y)\Big) = & \int_{\R^2} \partial_{\sigma}^l\Big(\log \mathcal L_{jk}(\sigma, y)\Big) \mathcal L_{jk}(\Sigma^*_{jk}, y)\mathrm{d}\mu(y) \\
= &\ \phi_{00,jk}(\Sigma^*_{jk})\partial_{\sigma}^l \log\phi_{00,jk}(\sigma) \\
& + \int_{a_k}^{b_k} \phi_{01,jk}(\Sigma^*_{jk}, y) \partial_{\sigma}^l\log \phi_{01,jk}(\sigma, y) \mathrm{d}y \\
& + \int_{a_j}^{b_j} \phi_{10,jk}(\Sigma^*_{jk}, x) \partial_{\sigma}^l\log \phi_{10,jk}(\sigma, x) \mathrm{d}x \\
& + \iint_{[a_j,b_j]\times[a_k,b_k]} \phi_{11,jk}(\Sigma^*_{jk}, x,y) \partial_{\sigma}^l\log \phi_{11,jk}(\sigma, x,y) \mathrm{d}x\mathrm{d}y.
\end{align*}
To show the equality $\partial_{\sigma}^l\E_{\Sigma^*_{jk}}\Big(\log \mathcal L_{jk}(\sigma, Y)\Big) = \E_{\Sigma^*_{jk}}\Big(\partial_{\sigma}^l \log \mathcal L_{jk}(\sigma, Y)\Big)$, we show the equality for each of the four terms

\vspace{0.2cm}
$\bullet$ For $a = 0, b = 0$: it is obvious. 

For the three remaining terms, we use Lebesgue theorem for differentiability of integrals with parameters.

\vspace{0.2cm}
$\bullet$ For $a = 0, b = 1$ (and $a = 1, b = 0$): 
$$\log \phi_{01,jk}(\sigma, y) = -\dfrac{y^2}{2} - \log\sqrt{2\pi} + \log \bigg[1 - F\Big(\dfrac{b_j-\sigma y}{\sqrt{1-\sigma^2}}\Big) + F\Big(\dfrac{a_j-\sigma y}{\sqrt{1-\sigma^2}} \Big) \bigg].$$ 
The function $\log \phi_{01,jk}$ is $C^\infty$ on the compact $[-1+\delta, 1-\delta]\times[a_k,b_k]$. Therefore, for all $\sigma \in [-1 + \delta, 1 - \delta]$ and $y \in [a_k, b_k]$, $\Big|\partial^l_{\sigma} \log \phi_{01,jk}(\sigma, y)\Big|$ is upper bounded by a constant 
which is integrable on the compact $[a_k, b_k]$ (with regard to the density $y \mapsto \phi_{01,jk}(\Sigma^*_{jk}, y)$).

\vspace{0.2cm}
$\bullet$ For $a = 1, b = 1$:
$$\log \phi_{11,jk}(\sigma, x, y) = - \log(2\pi) -\dfrac{1}{2}\log(1-\sigma^2)- \dfrac{x^2 - 2\sigma x y+ y^2}{2(1-\sigma^2)}.$$
Analogously, the function $\log \phi_{11,jk}$ is $C^\infty$ on the compact $[-1+\delta, 1-\delta]\times[a_j,b_j]\times[a_k,b_k]$. Thus, for all $\sigma \in [-1 + \delta, 1 - \delta]$, $x \in [a_j, b_j]$ and $y \in [a_k, b_k]$, $\Big|\partial^l_{\sigma} \log \phi_{11,jk}(\sigma, x, y)\Big|$ is upper bounded by a constant 
which is integrable on the compact $[a_j, b_j]\times[a_k, b_k]$ (with regard to the density $(x,y) \mapsto \phi_{11,jk}(\Sigma^*_{jk}, x,y)$). 
\end{enumerate}

\end{proof}

\section*{Acknowledgements}
We would like to thank Stéphane Robin and Stéphane Chrétien for their expert advice and comments.

\end{document}